\documentclass[12pt]{article}
\usepackage{amsfonts}
\usepackage[latin1]{inputenc}
\usepackage{amsmath,amssymb,amsfonts}
\usepackage{color}
\usepackage{hyperref}
\usepackage{enumerate}
\usepackage{graphicx}

\def\squarebox#1{\hbox to #1{\hfill\vbox to #1{\vfill}}}
\newcommand{\qed}{\hspace*{\fill}
\vbox{\hrule\hbox{\vrule\squarebox{.667em}\vrule}\hrule}\smallskip}
\newtheorem{teorema}{Theorem}[section]
\newtheorem{lema}[teorema]{Lemma}
\newtheorem{corolario}[teorema]{Corollary}
\newtheorem{proposicao}[teorema]{Proposition}

\newenvironment{profe}{\noindent {\bf Proof.}}{\hfill $\qed $ \newline}
\input{tcilatex}
\begin{document}

\title{A method to find generators of a semi-simple Lie group via the
topology of its flag manifolds}
\author{Ariane Luzia dos Santos\thanks{
Address: FCL - Unesp, Departamento de Ciências da Educação. Rodovia
Araraquara-Jaú, km 1, 14.800-901 Araraquara, São Paulo, Brasil. e-mail:
ariane@fclar.unesp.br} \and Luiz A. B. San Martin\thanks{
Supported by CNPq grant n$^{\mathrm{o}}$ 303755/2009-1 and FAPESP grant n$^{ 
\mathrm{o}}$ 07/06896-5}\thanks{
Address: Imecc - Unicamp, Departamento de Matemática. Rua Sérgio Buarque de
Holanda, 651, Cidade Universitária Zeferino Vaz, 13083-859 Campinas, São
Paulo, Brasil. e-mail: smartin@ime.unicamp.br}}
\date{}
\maketitle

\begin{abstract}
In this paper we continue to develop the topological method started in
Santos-San Martin \cite{ariasm} to get semigroup generators of semi-simple
Lie groups. Consider a subset $\Gamma \subset G$ that contains a semi-simple
subgroup $G_{1}$ of $G$. Then $\Gamma $ generates $G$ if $\mathrm{Ad}\left(
\Gamma \right) $ generates a Zariski dense subgroup of the algebraic group $%
\mathrm{Ad}\left( G\right) $. The proof is reduced to check that some
specific closed orbits of $G_{1}$ in the flag manifolds of $G$ are not
trivial in the sense of algebraic topology. Here, we consider three
different cases of semi-simple Lie groups $G$ and subgroups $G_{1}\subset G$.
\end{abstract}

\noindent%

\noindent%

\section{Introduction}

In this paper we continue to develop the topological method started in
Santos-San Martin \cite{ariasm} to get semigroup generators of semi-simple
Lie groups.

The method is based on the notion of flag type of a semigroup that arouse
from the results of \cite{smroconex}, \cite{sminv}, \cite{smsurv}, \cite{SM}%
, \cite{SMT}, \cite{smord}, \cite{smmax} and \cite{smsanhom}. By these
results if $G$ is a connected noncompact semi-simple Lie group with finite
center and $S\subset G$ is a proper semigroup with nonempty interior, then
the flag type of $S$ allows to select a flag manifold $\mathbb{F}_{\Theta }$
of $G$ that contains a subset $C_{\Theta }$ (the so called invariant control
set) which is invariant by $S$ and is contained in an open Bruhat cell $%
\sigma _{\Theta }$ of $\mathbb{F}_{\Theta }$. Since $\sigma _{\Theta }$ is
diffeomorphic to an Euclidean space $\mathbb{R}^{N}$ it follows that $%
C_{\Theta }$ is contractible in $\mathbb{F}_{\Theta }$.

Hence if one can show that a subset $\Gamma \subset G$ does not leave
invariant a contractible subset on any flag manifold of $G$, then it is
possible to conclude that $\Gamma $ generates $G$ if the semigroup generated
by $\Gamma $ has nonempty interior. Actually thanks to a result by Abels 
\cite{abel} this last condition can be changed by asking that $\mathrm{Ad}%
\left( \Gamma \right) $ generates a Zariski dense subgroup of the algebraic
group $\mathrm{Ad}\left( G\right) $.

The problem of semigroup generation of groups has several motivations. One
of them comes from control theory where the controllability problem is
translated into the semigroup generation problem. Controllability results on
semi-simple Lie groups where obtained in Jurdjevic-Kupka \cite{JK, JK1},
Gauthier-Kupka-Sallet \cite{GKS} and El Assoudi-Gauthier-Kupka \cite{EGK}.
In another direction we mention the 1.5 generation problem studied in Abels 
\cite{abel}, Abels-Vinberg \cite{abelsV} and references therein, which
consists in finding pairs of generators starting from one element of the
pair. 

Here, as in our previous paper \cite{ariasm}, we take as generator a subset $%
\Gamma \subset G$ that contains a semi-simple subgroup $G_{1}$ of $G$. In
this setting the problem is reduced to check that some specific closed
orbits of $G_{1}$ in the flag manifolds of $G$ are not trivial in the sense
of algebraic topology (see Proposition \ref{propSemiSimpleGeral} below).

In \cite{ariasm} we pursued this approach with an eye in the controllability
results mentioned above. Thus in \cite{ariasm} the group $G$ is a connected
complex simple Lie group and $G_{1}$ is a subgroup $G\left( \alpha \right) $
with Lie algebra isomorphic to $\mathfrak{sl}\left( 2,\mathbb{C}\right) $
generated by the root spaces $\mathfrak{g}_{\pm \alpha }$ associated to the
roots $\pm \alpha $ of the complex Lie algebra $\mathfrak{g}$ of $G$. In
that case the method was successfully applied because the relevant orbits of 
$G_{1}=G\left( \alpha \right) $ are flag manifolds of $\mathfrak{sl}\left( 2,%
\mathbb{C}\right) $ so that diffeomorphic to $S^{2}$. We used De Rham
cohomology $H^{2}\left( \mathbb{F}_{\Theta },\mathbb{R}\right) $ to prove
that these orbits are not homotopic to a point for any flag manifold $%
\mathbb{F}_{\Theta }$.

In this paper we consider three different cases of semi-simple Lie groups $G$
and subgroups $G_{1}\subset G$.

In the first one we take $G$ whose Lie algebra $\mathfrak{g}$ is the split
real form of a complex simple Lie algebra and $G_{1}=G\left( \alpha \right) $
is a subgroup generated by root spaces analogous to the complex case. In the
real case we are led to check whether some closed curves are null homotopic
and hence work with the fundamental groups of the flag manifolds. Contrary
to the complex Lie algebras there are only a few cases where the relevant
orbits of $G_{1}=G\left( \alpha \right) $ are null homotopic. Namely when $%
\mathfrak{g}=\mathfrak{sl}\left( n,\mathbb{R}\right) $, $\mathfrak{g}=%
\mathfrak{sp}\left( n,\mathbb{R}\right) $ and $\alpha $ is a long root and
when $\mathfrak{g}$ is the split real form of $G_{2}$ and $\alpha $ is a
short root. In these cases the subgroups \ $G\left( \alpha \right) $ are not
contained in proper semigroups with interior points and a subset $\Gamma $
generates $G$ if $G\left( \alpha \right) \subset \Gamma $ (any root $\alpha $%
) and the group generated by $\Gamma $ is Zariski dense (see Theorem \ref%
{teoreal} below). In the remaining cases our method breaks down. We give an
example of a proper semigroup with nonempty interior in $\mathrm{Sp}\left( n,%
\mathbb{R}\right) $ that contains $G\left( \alpha \right) $ for several
short roots $\alpha $, showing that the result is indeed not true in this
case.

In another direction we take an irreducible finite dimensional
representation $\rho _{n}:\mathfrak{sl}\left( 2,\mathbb{C}\right)
\rightarrow \mathfrak{sl}\left( n,\mathbb{C}\right) $, $n\geq 2$, and get a
\linebreak subgroup $G_{1}=\langle \exp \rho _{n}\left( \mathfrak{sl}\left(
2,\mathbb{C}\right) \right) \rangle \subset \mathrm{Sl}\left( n,\mathbb{C}%
\right) $. We prove that $G_{1}$ is not contained in any proper semigroup
with nonempty interior of $\mathrm{Sl}\left( n,\mathbb{C}\right) $. It
follows that a subset $\Gamma $ is a semigroup generator of $\mathrm{Sl}%
\left( n,\mathbb{C}\right) $ if it contains $\langle \exp \rho _{n}\left( 
\mathfrak{sl}\left( 2,\mathbb{C}\right) \right) \rangle $ and the group
generated by $\Gamma $ is Zariski dense (see Theorem \ref{teoParaRepreSl2}).
The algebraic topological fact that permits the proof of this result is
Proposition \ref{propfiblinhas2}. It shows that the projective orbit $%
G_{1}\cdot \left[ v_{0}\right] $ of the highest weight space is not
contractible in the projective space of the vector space of the
representation. This noncontractibility follows from the fact that the

tautological bundle of the projective space restricts to a nontrivial bundle
on the orbit

In our third case we take a complex Lie group $G$ and a complex subgroup $%
G_{1}$ such that the Lie algebra $\mathfrak{g}_{1}$ of $G_{1}$ contains a
regular real element of the Lie algebra $\mathfrak{g}$ of $G$. Examples of
this case are the inclusion in $\mathfrak{sl}\left( n,\mathbb{C}\right) $ of
the classical Lie algebras $B_{l}=\mathfrak{so}\left( 2l+1,\mathbb{C}\right) 
$, $C_{l}=\mathfrak{sp}\left( 2l,\mathbb{C}\right) $ and $D_{l}=\mathfrak{so}%
\left( 2l,\mathbb{C}\right) $. In this case we prove that $G_{1}$ is not
contained in a proper semigroup with nonempty interior of $G$ and hence get
generators of $G$ in the same spirit of the other cases (see Theorem \ref%
{teoParaSemiComplex}). Here we exploit the same technique provided by
Proposition \ref{propfiblinhas2} by realizing the flag manifolds of $G$ as
projective orbits in spaces of representations. The proof is not much
different from the case of representations of $\mathfrak{sl}\left( 2,\mathbb{%
C}\right) $.

\section{Notation and background}

Let $G$ be a connected real semi-simple Lie group with finite center and Lie
algebra $\mathfrak{g}$. For $G$ and $\mathfrak{g}$ we use the following
notation.

\begin{itemize}
\item Let $\theta $ be a Cartan involution of $\mathfrak{g}$ and $\mathfrak{g%
}=\mathfrak{k}\oplus \mathfrak{s}$ the corresponding Cartan decomposition.

\item If $\mathfrak{a}\subset \mathfrak{s}$ is a maximal abelian subalgebra
its set of roots is denoted by $\Pi $. Let $\Pi ^{+}$ be a set of positive
roots with 
\begin{equation*}
\Sigma =\{\alpha _{1},\ldots ,\alpha _{l}\}\subset \Pi ^{+}
\end{equation*}
corresponding simple system of roots. We have $\Pi =\Pi ^{+}\dot{\cup}(-\Pi
^{+})$ and any $\alpha \in \Pi ^{+}$ is a linear combination $\alpha
=n_{1}\alpha _{1}+\cdots +n_{l}\alpha _{l}$ and $n_{i}\geq 0$, $i=1,\ldots
,l $ are integers. The support of $\alpha $, $\mathrm{supp}\alpha $, is the
subset of $\Sigma $ where $n_{i}>0$.

\item The Cartan-Killing form of $\mathfrak{g}$ is denoted by $\langle \cdot
,\cdot \rangle $. If $\alpha \in \mathfrak{a}^{\ast }$ then $H_{\alpha }\in 
\mathfrak{a}$ is defined by $\alpha (\cdot )=\langle H_{\alpha },\cdot
\rangle $, and $\langle \alpha ,\beta \rangle =\langle H_{\alpha },H_{\beta
}\rangle $.

\item We write 
\begin{equation*}
\mathfrak{a}^{+}=\{H\in \mathfrak{a}:\forall \alpha \in \Pi ^{+},\,\alpha
(H)>0\}
\end{equation*}
for the Weyl chamber defined by $\Pi ^{+}$.

\item The root space of a root $\alpha $ is 
\begin{equation*}
\mathfrak{g}_{\alpha }=\{X\in \mathfrak{g}:\forall H\in \mathfrak{a}%
,\,[H,X]=\alpha (H)X\}.
\end{equation*}
If $\mathfrak{g}$ is a split real form of a complex simple Lie algebra then $%
\dim _{\mathbb{R}}\mathfrak{g}_{\alpha }=1$.

\item For a root $\alpha $, $\mathfrak{g}(\alpha )$ is the subalgebra
generated by $\mathfrak{g}_{\alpha }$ and $\mathfrak{g}_{-\alpha }$. We have 
\begin{equation*}
\mathfrak{g}(\alpha )=\mathrm{span}_{\mathbb{R}}\{H_{\alpha }\}\oplus 
\mathfrak{g}_{\alpha }\oplus \mathfrak{g}_{-\alpha }\approx \mathfrak{sl}(2,%
\mathbb{R}).
\end{equation*}
We let $G(\alpha )$ be the connected Lie subgroup of $G$ with Lie algebra $%
\mathfrak{g}(\alpha )$.

\item Let $K=\langle \exp \mathfrak{k}\rangle $ be the maximal compact
subgroup of $G$ with Lie algebra $\mathfrak{k}$.

\item $\mathcal{W}$ is the Weyl group. Either $\mathcal{W}$ is the group
generated by the reflections $r_{\alpha }$, $\alpha \in \Pi $, $r_{\alpha
}(\beta )=\beta -\frac{2\langle \alpha ,\beta \rangle }{\langle \alpha
,\alpha \rangle }\alpha $, or $\mathcal{W}=M^{\ast }/M$ where $M^{\ast }=%
\mathrm{Norm}_{K}(\mathfrak{a})$ is the normalizer of $\mathfrak{a}$ in $K$
and $M=\{g\in K:\mathrm{Ad}(g)H=H\,\mbox{for all}\,H\in \mathfrak{h}\}$ is
the centralizer of $\mathfrak{a}$ in $K$.

\item $\mathfrak{n}^{+}=\sum_{\alpha \in \Pi ^{+}}\mathfrak{g}_{\alpha }$
and $\mathfrak{n}^{-}=\sum_{\alpha \in \Pi ^{+}}\mathfrak{g}_{-\alpha }$.

\item Given the data $\mathfrak{a}$ and $\Pi ^{+}$ (or $\Sigma $) there
exists the minimal parabolic subalgebra $\mathfrak{p}=\mathfrak{m}\oplus 
\mathfrak{a}\oplus \mathfrak{n}^{+}$. A subset $\Theta \subset \Sigma $
defines the standard parabolic subalgebra by 
\begin{equation*}
\mathfrak{p}_{\Theta }=\mathfrak{p}+\sum_{\alpha \in \langle \Theta \rangle }%
\mathfrak{g}_{\alpha }
\end{equation*}
where $\langle \Theta \rangle =\{\alpha \in \Pi :$ $\mathrm{supp}\alpha
\subset \Theta $ or $\mathrm{supp}(-\alpha )\subset \Theta \}$ is the set of
roots spanned by $\Theta $. When $\Theta =\emptyset $ we have $\mathfrak{p}%
_{\emptyset }=\mathfrak{p}$.

\item For $\Theta \subset \Sigma $, $P_{\Theta }$ is the parabolic subgroup
with Lie algebra $\mathfrak{p}_{\Theta }$: 
\begin{equation*}
P_{\Theta }=\mathrm{Norm}_{G}(\mathfrak{p}_{\Theta })=\{g\in G:\mathrm{Ad}(g)%
\mathfrak{p}_{\Theta }\subset \mathfrak{p}_{\Theta }\}.
\end{equation*}

\item The flag manifold $\mathbb{F}_{\Theta }=G/P_{\Theta }$ does not depend
on the specific group $G$ with Lie algebra $\mathfrak{g}$. The origin of $%
G/P_{\Theta }$, the coset $1\cdot P_{\Theta }$, is denoted by $b_{\Theta }$.
\end{itemize}

Now let $S\subset G$ be a subsemigroup with $\mathrm{int}S\neq \emptyset $.
We recall some results of \cite{sminv, smsurv, SMT} that are on the basis of
our topological approach to get generators of $G$.

We let $S$ act on a flag manifold $\mathbb{F}_{\Theta }$ by restricting the
action of $G$. An invariant control set for $S$ in $\mathbb{F}_{\Theta }$ is
a subset $C\subset \mathbb{F}_{\Theta }$ such that $\mathrm{cl}(Sx)=C$ for
every $x\in C$, where $Sx=\{gx\in \mathbb{F}_{\Theta }:g\in S\}$. Since $%
\mathrm{int}S\neq \emptyset $ such a set is closed, has nonempty interior
and is in fact invariant, that is, $gx\in C$ if $g\in S$ and $x\in C$.

\begin{lema}
\textrm{{(\cite[Theorem 3.1]{sminv})}} \label{lemuniqueics}In any flag
manifold $F_{\Theta }$ there is a unique invariant control set for $S$,
denoted by $C_{\Theta }$.
\end{lema}

To state the geometric property of $C_{\Theta }$ to be used later we discuss
the dynamics of the vector fields $\widetilde{H}$ on a flag manifold $%
\mathbb{F}_{\Theta }$ whose flow is $e^{tH}$, with $H$ in the closure $%
\mathrm{cl}\mathfrak{a}^{+}$ of the Weyl chamber $\mathfrak{a}^{+}$. It is
known that $\widetilde{H}$ is a gradient vector field with respect to some
Riemmannian metric on $\mathbb{F}_{\Theta }$; cf. \cite[Proposition 3.3 (ii)]%
{dkv} and \cite{fepase}.

Hence the orbits of $\widetilde{H}$ are either fixed points of $e^{tH}$ or
trajectories flowing between fixed points of $H$. Moreover, $\widetilde{H}$
has a unique attractor fixed point set, say $\mathrm{att}_{\Theta }(H)$,
that has an open and dense stable manifold $\sigma _{\Theta }(H)$; cf. \cite%
{dkv, fepase}. This means that if $x\in \sigma _{\Theta }(H)$ then its $%
\omega $-limit set $\omega (x)$ is contained in $\mathrm{att}_{\Theta }(H)$.
This attractor has the following algebraic expressions 
\begin{equation*}
\mathrm{att}_{\Theta }(H)=Z_{H}\cdot b_{\Theta }=K_{H}\cdot b_{\Theta },
\end{equation*}
cf. \cite[Corollary 3.5]{dkv} and \cite{fepase}. Here $Z_{H}=\{g\in G:%
\mathrm{Ad}(g)H=H\}$ is the centralizer of $H$ in $G$ and $K_{H}=Z_{H}\cap K$
is the centralizer in $K$. We note that $\mathrm{att}_{\Theta }(H)$ is a
connected manifold because $Z_{H}=M(Z_{H})_{0}=(Z_{H})_{0}M$ where $%
(Z_{H})_{0}$ is the identity component of $Z_{H}$ and $M$ is the centralizer
of $\mathfrak{a}$ in $K$; see \cite[Lemma 1.2.4.5]{W}. Hence $Z_{H}\cdot
b_{\Theta }=(Z_{H})_{0}\cdot b_{\Theta }$, since $M\cdot b_{\Theta
}=b_{\Theta }$; cf. \cite[Theorem 7.101]{K}.

The stable set $\sigma _{\Theta }(H)$ is also described algebraically by 
\begin{equation*}
\sigma _{\Theta }(H)=N_{H}^{-}Z_{H}\cdot b_{\Theta }
\end{equation*}
where 
\begin{equation*}
N_{H}^{-}=\exp \mathfrak{n}_{H}^{-}\quad \mathrm{and}\quad \mathfrak{n}%
_{H}^{-}=\sum_{\gamma (H)<0}\mathfrak{g}_{\gamma },
\end{equation*}
cf. \cite[Corollary 3.5]{dkv}. In particular if $H$ is regular, that is, $%
H\in \mathfrak{a}$ and $\alpha (H)>0$ for $\alpha \in \Pi ^{+}$, then $%
Z_{H}=MA$, which fixes $b_{\Theta }$. Hence 
\begin{equation*}
\mathrm{att}_{\Theta }(H)=Z_{H}\cdot b_{\Theta }=\{b_{\Theta }\}\qquad H\in 
\mathfrak{a}.
\end{equation*}
Actually, in the regular case the fixed points are isolated because $%
\widetilde{H}$ is the gradient of a Morse function; cf. \cite{dkv, fepase}.
Also, $\mathfrak{n}_{H}^{-}=\mathfrak{n}^{-}$ (notation as above) and the
stable set is $N^{-}\cdot b_{\Theta }$ the open Bruhat cell.

The following statement is a well known result from the Bruhat decomposition
of the flag manifolds; cf. \cite{dkv, K, W}.

\begin{proposicao}
In any flag manifold $\mathbb{F}_{\Theta }$ the open Bruhat cell $N^{-}\cdot
b_{\Theta }$ is diffeomorphic to an Euclidean space $\mathbb{R}^{d}$. The
diffeomorphism is $X\in \mathfrak{n}_{\Theta }^{-}\mapsto e^{X}\cdot
b_{\Theta }$, where $\mathfrak{n}_{\Theta }^{-}=\sum \{\mathfrak{g}_{\alpha
}:\alpha <0$ and $\alpha \notin \langle \Theta \rangle \}$.
\end{proposicao}

Set $h=e^{H}$, $H\in \mathfrak{a}^{+}$. It follows from the gradient
property of $\widetilde{H}$ that $\lim_{n\rightarrow +\infty
}h^{n}x=b_{\Theta }$ for any $x\in N^{-}\cdot b_{\Theta }$.

Now, we say that $g\in G$ is regular real if it is a conjugate $g=aha^{-1}$
of $h=\exp H$, $H\in \mathfrak{a}^{+}$ with $a\in G$. Then we write $\sigma
_{\Theta }( g) =g\cdot \sigma _{\Theta}( H) $ and we call this the stable
set of $g$ in $\mathbb{F}_{\Theta }$. The reason for this name is clear: $%
g^{n}=(aha^{-1}) ^{n}=ah^{n}a^{-1}$ and hence $g^{n}x\rightarrow gb_{\Theta
} $ if $x\in \sigma _{\Theta }( g) $.

The following lemma was used in \cite{sminv} to prove the above Lemma \ref%
{lemuniqueics}.

\begin{lema}
\textrm{{(\cite[Lemma 3.2]{sminv})}} There exists a regular real $g\in%
\mathrm{int}S$.
\end{lema}

Now we can state the following result of \cite{SMT} which is basic to our
approach.

\begin{teorema}
\label{teoregre}Suppose that $S\neq G$. Then there exists a flag manifold $%
\mathbb{F}_{\Theta }$ such that the invariant control set $C_{\Theta
}\subset \sigma _{\Theta }(g)$ for every regular real $g\in \mathrm{int}S$.
\end{teorema}

\begin{corolario}
\label{coreuclid}If $S\neq G$ then there exists a flag manifold $\mathbb{F}%
_{\Theta }$ such that for every flag manifold $\mathbb{F}_{\Theta _{1}}$
such that $\Theta \subset \Theta _{1}$ the invariant control set $C_{\Theta
_{1}}$ in $\mathbb{F}_{\Theta _{1}}$ is contained in a subset $\mathcal{E}%
_{\Theta _{1}}$ diffeomorphic to an Euclidean space.
\end{corolario}

\begin{profe}
If $\Theta _{1}\subset \Theta $ then the canonical projection $\pi :\mathbb{F%
}_{\Theta }\rightarrow \mathbb{F}_{\Theta _{1}}$ is equivariant under the
actions of $G$. This implies that the open Bruhat cells are projected onto
open cells\ and $\pi \left( C_{\Theta }\right) =C_{\Theta _{1}}$. Hence $%
C_{\Theta _{1}}$ is contained in an open cell $\mathcal{E}_{\Theta _{1}}$ if
this happens to $C_{\Theta }$.
\end{profe}

In particular if $\Theta _{1}=\Sigma \setminus \{\alpha \}$ contains $\Theta 
$ if $\alpha \in \Sigma \setminus \Theta $ so in the minimal flag manifold $%
\mathbb{F}_{\Sigma \setminus \{\alpha \}}$ the invariant control set is
contained in open cells.

\begin{corolario}
\label{coreuclidminimal}If $S\neq G$ then there exists a minimal flag
manifold $\mathbb{F}_{\Theta }$ such that $C_{\Theta }$ is contained in a
subset $\mathcal{E}_{\Theta }$ diffeomorphic to an Euclidean space.
\end{corolario}

\vspace{12pt}%

\noindent \textbf{Remark:} It can be proved that there exists a minimal $%
\Theta \left( S\right) $ satisfying the condition of Theorem \ref{teoregre}.
This $\Theta \left( S\right) $ (or rather the flag manifold $\mathbb{F}%
_{\Theta \left( S\right) }$) is called the flag type or parabolic type of $S$%
. Several properties of $S$ are derived from this flag type, e.g. the
homotopy type of $S$ as in \cite{smsanhom} or the connected components of $S$
as in \cite{smroconex}.

\section{Semi-simple subgroups: Set up\label{secSemiGeral}}

Let $\mathfrak{g}$ be a noncompact semi-simple Lie algebra and $G$ be a
connected Lie group with Lie algebra $\mathfrak{g}$ and finite center. In
this section we take a semi-simple subalgebra $\mathfrak{g}_{1}\subset 
\mathfrak{g}$ and the corresponding connected subgroup $G_{1}=\langle \exp 
\mathfrak{g}_{1}\rangle $. We ask whether there is a proper semigroup $%
S\subset G$ with $\mathrm{int}S\neq \emptyset $ such that $G_{1}\subset S$.

It is well known that there exist compatible Cartan decompositions $%
\mathfrak{g}_{1}=\mathfrak{k}_{1}\oplus \mathfrak{s}_{1}$ and $\mathfrak{g}=%
\mathfrak{k}\oplus \mathfrak{s}$ such that $\mathfrak{k}_{1}\subset 
\mathfrak{k}$ and $\mathfrak{s}_{1}\subset \mathfrak{s}$ (see Warner \cite[%
Lemma 1.1.5.5]{W}). If $\mathfrak{a}_{1}\subset \mathfrak{s}_{1}$ is maximal
abelian then there exists a maximal abelian $\mathfrak{a}\subset \mathfrak{s}
$ with $\mathfrak{a}_{1}\subset \mathfrak{a}$. Denote by $\Pi _{1}$ the
roots of $\left( \mathfrak{g}_{1},\mathfrak{a}_{1}\right) $ and by $\Pi $
the roots of $\left( \mathfrak{g},\mathfrak{a}\right) $. Any $\alpha _{1}\in
\Pi _{1}$ is the restriction to $\mathfrak{a}_{1}$ of some $\alpha \in \Pi $%
. Take $H_{1}\in \mathfrak{a}_{1}$ regular (in $\mathfrak{g}_{1}$). Then
there exists a Weyl chamber $\mathfrak{a}^{+}\subset \mathfrak{a}$ with $%
H_{1}\in \mathrm{cl}\mathfrak{a}^{+}$. If $\alpha \in \Pi $ is positive
w.r.t $\mathfrak{a}^{+}$ then $\alpha \left( H_{1}\right) \geq 0$ so that we
get compatible Iwasawa decompositions $\mathfrak{g}_{1}=\mathfrak{k}%
_{1}\oplus \mathfrak{a}_{1}\oplus \mathfrak{n}_{1}$ and $\mathfrak{g}=%
\mathfrak{k}\oplus \mathfrak{a}\oplus \mathfrak{n}$ with $\mathfrak{k}%
_{1}\subset \mathfrak{k}$, $\mathfrak{a}_{1}\subset \mathfrak{a}$ and $%
\mathfrak{n}_{1}\subset \mathfrak{n}$ where 
\begin{equation*}
\mathfrak{n}_{1}=\sum_{\alpha _{1}\in \Pi _{1}^{+}}\left( \mathfrak{g}%
_{1}\right) _{\alpha _{1}}\, ,\quad \mathfrak{n}=\sum_{\alpha \in \Pi ^{+}}%
\mathfrak{g}_{\alpha }
\end{equation*}
and $\Pi _{1}^{+}$ is the set of roots of $\left( \mathfrak{g}_{1},\mathfrak{%
a}_{1}\right) $ that are positive on $H_{1}$ and $\Pi ^{+}$ the roots
positive on $\mathfrak{a}^{+}$.

In the sequel we keep fixed these compatible Iwasawa decompositions. Let $%
\Sigma $ be the simple system of roots in $\Pi ^{+}$. Then the standard
parabolic subalgebras $\mathfrak{p}_{\Theta }\subset \mathfrak{g}$ and
subgroups $P_{\Theta }\subset G$ are built from subsets $\Theta
\subset\Sigma $. For the corresponding flag manifolds $\mathbb{F}%
_{\Theta}=G/P_{\Theta }$ we write $b_{\Theta }=1\cdot P_{\Theta }$ for their
origins.

By the construction of the compatible Iwasawa decomposition from the choice
of $H_{1}\in \mathfrak{a}_{1}\cap \mathrm{cl}\mathfrak{a}^{+}$ we have that $%
b_{\Theta }$ belongs to the attractor fixed point set $\mathrm{att}%
_{\Theta}(H_{1})=Z_{H_{1}}\cdot b_{\Theta }=K_{H_{1}}\cdot b_{\Theta }$ of
the one-parameter semigroup $e^{tH_{1}}$, $t\geq 0$. The corresponding
stable set $\sigma _{\Theta }(H_{1})$ is open and dense in $\mathbb{F}%
_{\Theta }$.

Now let $S$ be a semigroup with $\mathrm{int}S\neq \emptyset $ such that $%
G_{1}\subset S$. Denote by $C_{\Theta }$ the unique $S$-invariant control
set in $\mathbb{F}_{\Theta }$. The set $C_{\Theta }$ is compact and has
nonempty interior. Hence $C_{\Theta }\cap \sigma _{\Theta }(H_{\alpha
})\neq\emptyset $. If $x\in C_{\Theta }\cap \sigma _{\Theta }(H_{\alpha })$
then $y=\lim_{t\rightarrow +\infty }e^{tH_{\alpha }}\cdot x$ belongs to $%
C_{\Theta}\cap \mathrm{att}_{\Theta }(H_{\alpha })$, because $C_{\Theta }$
is closed. Hence we get the following lemma.

\begin{lema}
\label{lemintericsatrac}If $G_{1}\subset S$ then $C_{\Theta }\cap \mathrm{att%
}_{\Theta }(H_{1})\neq \emptyset $.
\end{lema}

Now the idea is to look at the topology of the orbits $G_{1}\cdot y$ with $%
y\in \mathrm{att}_{\Theta }(H_{1})$. Clearly if $y\in C_{\Theta }$ and $%
G_{1}\subset S$ then $G_{1}\cdot y\subset C_{\Theta }$. On the other hand,
by Theorem \ref{teoregre} and its corollaries, there are flag manifolds
where $C_{\Theta }$ is contained in a contractible Euclidean subset $%
\mathcal{E}_{\Theta }$ if $S$ is proper. Hence if we achieve to prove that
none of the orbits $G_{1}\cdot y$ with $y\in \mathrm{att}_{\Theta }(H_{1})$
are contractible then we can conclude that $G_{1}$ is not contained in a
proper semigroup with nonempty interior. In principle this
noncontractibility property must by checked on every flag manifold but by
Corollary \ref{coreuclidminimal} it is enough to look at the minimal ones.

These arguments can be used to get semigrouop generators of $G$. In fact, if 
$\Gamma $ is a subset that contains $G_{1}$ and generates a semigroup $S$
with nonempty interior then $S=G$ provided we have noncontractibility of the
orbits $G_{1}\cdot y$ through the attractor fixed point set.

Actually, thanks to a result by Abels \cite{abel} it is enough to assume
that the group generated by $\Gamma $ is Zariski dense in the following
sense: The group $\mathrm{Ad}\left( G\right) $ is algebraic and hence
endowed with the Zariski topology. We say that $B\subset G$ is Zariski dense
in case $\mathrm{Ad}\left( B\right) $ is dense in $\mathrm{Ad}\left(
G\right) $ with respect to the Zariski topology. With this terminology it is
proved in \cite{abel}, as a consequence of Corollary 5, that the semigroup $%
S $ generated by $\Gamma $ has non empty interior provided i) the group
generated by $\Gamma $ is Zariski dense and ii) $S$ contains a non-constant
smooth curve.

These comments yield the following fact that reduces the problem of finding
semigroup generators to the topology of orbits of $G_{1}$.

\begin{proposicao}
\label{prop:zariski} \label{propSemiSimpleGeral}Let $\mathfrak{g}_{1}\subset 
\mathfrak{g}$ be a semi-simple Lie subalgebra. Choose compatible Iwasawa
decompositions $\mathfrak{g}_{1}=\mathfrak{k}_{1}\oplus \mathfrak{a}%
_{1}\oplus \mathfrak{n}_{1}\subset \mathfrak{g}=\mathfrak{k}\oplus \mathfrak{%
a}\oplus \mathfrak{n}$, so that $\mathrm{cl}\mathfrak{a}^{+}$ contains a
regular element $H_{1}\in \mathfrak{a}_{1}$. Let $G_{1}=\langle \exp 
\mathfrak{g}_{1}\rangle $ and suppose that for every (minimal) flag manifold 
$\mathbb{F}_{\Theta }$ the orbits $G_{1}\cdot y$ through the attractor fixed
point set $\mathrm{att}_{\Theta }(H_{1})$ are not contractible in $\mathbb{F}%
_{\Theta }$. Then a subset $\Gamma \subset G$ generates $G$ as a semigroup
if $G_{1}\subset \Gamma $ and the subgroup generated by $\Gamma $ is Zariski
dense.
\end{proposicao}

In the special case when $H_{1}\in \mathfrak{a}^{+}$, that is, $\mathfrak{g}%
_{1}$ contains a regular real element of $\mathfrak{g}$, the attractor fixed
point $\mathrm{att}_{\Theta }(H_{1})$ reduces to $b_{\Theta }$. In this case
we need to check contractibility only of the orbits $G_{1}\cdot b_{\Theta }$
through the origin. For later reference we record this fact.

\begin{corolario}
\label{corSemiSimpleGeral}With the notation of Proposition \ref{prop:zariski}
suppose that $H_{1}\in \mathfrak{a}^{+}$. Then the same result holds with
the assumption that $G_{1}\cdot b_{\Theta }$ is not contractible in any
(minimal) flag manifold.
\end{corolario}

\section{Split real forms and subgroups $G\left( \protect\alpha \right) $}

We assume in this section that $\mathfrak{g}$ is a split real form of a
complex simple Lie algebra and $G$ is a connected Lie group with Lie algebra 
$G$ and having finite center. Take a Cartan decomposition $\mathfrak{g}=%
\mathfrak{k}\oplus \mathfrak{s}$ and a maximal abelian $\mathfrak{a}\subset 
\mathfrak{s}$. If $\alpha $ is a root of the pair $\left( \mathfrak{g},%
\mathfrak{a}\right) $ then the subalgebra $\mathfrak{g}\left( \alpha \right) 
$ generated by the root spaces $\mathfrak{g}_{\pm \alpha }$ is isomorphic to 
$\mathfrak{sl}\left( 2,\mathbb{R}\right) $. Precisely, we choose $X_{\alpha
}\in \mathfrak{g}_{\alpha }$, $X_{-\alpha }=-\theta X_{\alpha }\in \mathfrak{%
g}_{-\alpha }$ such that $\langle X_{\alpha },X_{-\alpha }\rangle =1$. Then
the isomorphism is given by 
\begin{equation*}
X_{\alpha }\leftrightarrow \left( 
\begin{array}{cc}
0 & 1 \\ 
0 & 0%
\end{array}
\right) ,\qquad H_{\alpha }\leftrightarrow \left( 
\begin{array}{cc}
1 & 0 \\ 
0 & -1%
\end{array}
\right) ,\qquad X_{-\alpha }\leftrightarrow \left( 
\begin{array}{cc}
0 & 0 \\ 
1 & 0%
\end{array}
\right) .
\end{equation*}

In order to apply the ideas of Section \ref{secSemiGeral} we fix a Weyl
chamber $\mathfrak{a}^{+}$ such that $H_{\alpha }\in \mathrm{cl}\mathfrak{a}%
^{+}$. Then $\alpha $ is positive w.r.t. $\mathfrak{a}^{+}$ and we get
compatible Iwasawa decompositions $\mathfrak{g}\left( \alpha \right) =%
\mathfrak{k}_{H_{\alpha }}\oplus \langle H_{\alpha }\rangle \oplus \mathfrak{%
g}_{\alpha }\subset \mathfrak{g}=\mathfrak{k}\oplus \mathfrak{a}\oplus 
\mathfrak{n}$ where $\mathfrak{k}_{H_{\alpha }}$ is spanned by $%
A_{\alpha}=X_{\alpha }-X_{-\alpha }$. We put $G\left( \alpha \right)
=\langle \exp \mathfrak{g}\left( \alpha \right) \rangle $.

Our objective is to check whether the conditions of Proposition \ref%
{propSemiSimpleGeral} are satisfied by $G\left( \alpha \right) $. The
following result reduces the question to the orbit through the origin $%
b_{\Theta }$ of a flag manifold $\mathbb{F}_{\Theta }$. (Where the origins
are given by the standard parabolic subalgebras defined from the Weyl
chamber $\mathfrak{a}^{+}$.)

\begin{teorema}
\label{teoOrbitsAtt}In a flag manifold $\mathbb{F}_{\Theta }$ take $y\in 
\mathrm{att}_{\Theta }(H_{\alpha })$. Then the orbit $G(\alpha )\cdot y$ is
a circle $S^{1}$ homotopic to the orbit $G(\alpha )\cdot b_{\Theta }$.
\end{teorema}

The proof this theorem requires some lemmas. We start by looking at the
orbit $G(\alpha )\cdot b_{\Theta }$. By compatibility of the Iwasawa
decompositions it follows that the parabolic subalgebra $\mathfrak{p}%
_{\alpha }=\langle H_{\alpha }\rangle \oplus \mathfrak{g}_{\alpha }$ of $%
\mathfrak{g}\left( \alpha \right) $ is contained in the isotropy subalgebra
at $b_{\Theta }$ for any flag manifold $\mathbb{F}_{\Theta }$. To get the
inclusion of the parabolic subgroup of $G\left( \alpha \right) $ as well we
perform a somewhat standard computation on $\mathfrak{sl}(2,\mathbb{R})$ and 
$\mathfrak{sl}\left( 2,\mathbb{C}\right) $ (cf. \cite{K}, Chapter VII.5).

\begin{proposicao}
\label{propcentroreal}Let $G(\alpha )\subset G$ be the connected subgroup
with Lie algebra $\mathfrak{g}(\alpha )\approx \mathfrak{sl}(2,\mathbb{R})$
and $M$ the centralizer of $\mathfrak{a}$ in $K$. Then the center $Z\left(
G(\alpha )\right) $ of $G(\alpha )$ is contained in $M$.
\end{proposicao}

\begin{profe}
Let $\widetilde{G}=\widetilde{\mathrm{Sl}(2,\mathbb{R})}$ be the universal
covering and denote by $\widetilde{\exp }:\mathfrak{sl}(2,\mathbb{R}%
)\rightarrow \widetilde{\mathrm{Sl}(2,\mathbb{R})}$ the exponential map.
Take the basis 
\begin{equation*}
A=\left( 
\begin{array}{cc}
0 & -1 \\ 
1 & 0%
\end{array}
\right) ,\qquad H=\left( 
\begin{array}{cc}
1 & 0 \\ 
0 & -1%
\end{array}
\right) ,\qquad S=\left( 
\begin{array}{cc}
0 & 1 \\ 
1 & 0%
\end{array}
\right)
\end{equation*}
of $\mathfrak{sl}(2,\mathbb{R})$. The center $Z\left( \widetilde{\mathrm{Sl}
(2,\mathbb{R})}\right) $ of $\widetilde{\mathrm{Sl}(2\mathbb{R})}$ is the
kernel of its adjoint representation $\widetilde{\mathrm{Ad}}$ which is
explicitly given by 
\begin{equation*}
Z\left( \widetilde{\mathrm{Sl}(2,\mathbb{R})}\right) =\{\widetilde{\exp }
(k\pi A):k\in \mathbb{Z}\}\approx \mathbb{Z},
\end{equation*}
since the center is contained in a one-parameter group $\widetilde{\exp }%
(tA) $.

Now $\mathfrak{g}(\alpha )\approx \mathfrak{sl}(2,\mathbb{R})$ with the
isomorphism given by $H\leftrightarrow H_{\alpha }^{\vee }=\frac{2H_{\alpha }%
}{\langle \alpha ,\alpha \rangle }$, $A\leftrightarrow A_{\alpha }=X_{\alpha
}-X_{-\alpha }$ and $S\leftrightarrow S_{\alpha }=X_{\alpha }+X_{-\alpha }$,
where $X_{\pm \alpha }\in \mathfrak{g}_{\pm \alpha }$ and $\langle X_{\alpha
},X_{-\alpha }\rangle =1$.

Suppose first that $G=\mathrm{Aut}_{0}\mathfrak{g}$ is the adjoint group.
Then $G(\alpha )=\widetilde{\mathrm{Sl}(2,\mathbb{R})}/D$ with $D\subset
Z\left( \widetilde{\mathrm{Sl}(2,\mathbb{R})}\right) $ given by 
\begin{equation*}
D=\{\widetilde{\exp }(kn_{0}\pi A):k\in \mathbb{Z}\}\approx n_{0}\mathbb{Z}
\end{equation*}
where 
\begin{equation*}
n_{0}=\min \{n>0:e^{n\pi \mathrm{ad}\left( A_{\alpha }\right) }=\mathrm{id}
\}.
\end{equation*}
It follows that $Z\left( G(\alpha )\right) =Z\left( \widetilde{\mathrm{Sl}(2,%
\mathbb{R})}\right) /D\approx \mathbb{Z}_{n_{0}}$ is generated by $e^{\pi 
\mathrm{ad}\left( A_{\alpha }\right) }$. Complexifying and doing
computations in $\mathrm{Sl}(2,\mathbb{C})$, it turns out that $e^{\pi 
\mathrm{ad}(A_{\alpha })}=e^{\pi \mathrm{ad}(iH_{\alpha }^{\vee })}$. This
last term belongs to $M$ showing that $Z\left( G(\alpha )\right) \subset M$
when $G=\mathrm{Aut}_{0}\mathfrak{g}$. For a general $G$ the same result is
obtained by taking adjoints.
\end{profe}

The next lemma will ensure that the isotropy subalgebra at $b_{\Theta }$ for
the action of $G\left( \alpha \right) $ is exactly $\mathfrak{p}_{\alpha
}=\langle H_{\alpha }\rangle \oplus \mathfrak{g}_{\alpha }$.

\begin{lema}
\label{lemsuport}If $\alpha $ is a root with $H_{\alpha }\in \mathrm{cl}%
\mathfrak{a}^{+}$, then $\mathrm{supp}\alpha =\Sigma $.
\end{lema}

\begin{profe}
See \cite[Proposition 3.3]{ariasm}.
\end{profe}

Now we can describe the orbits through the origins of the flag manifolds.

\begin{lema}
\label{lemesse1}Let $b_{\Theta }$ be the origin of a flag manifold $\mathbb{F%
}_{\Theta }$ and $\beta $ a positive root. Then $G(\beta )\cdot b_{\Theta }$
is either a circle $S^{1}$ or it reduces to a point. If $\beta \notin
\langle \Theta \rangle $ then $\dim \,G(\beta )\cdot b_{\Theta }=1$. In
particular, $\dim G(\alpha )\cdot b_{\Theta }=1$ if $H_{\alpha }\in \mathrm{%
cl}\mathfrak{a}^{+}$.
\end{lema}

\begin{profe}
Let $\mathfrak{g}\left( \beta \right) _{b_{\Theta }}$ and $G\left( \beta
\right) _{b_{\Theta }}$ be the isotropy subalgebra and subgroup,
respectively, at $b_{\Theta }$ for the action of $G\left( \beta \right) $.
The subalgebra $\mathfrak{g}\left( \beta \right) _{b_{\Theta }}$ contains
the parabolic subalgebra of $\mathfrak{g}\left( \beta \right) $ given by $%
\mathfrak{p}_{\beta }=\mathrm{span}\{H_{\beta }\}\oplus \mathfrak{g}_{\beta
}\subset \mathfrak{p}_{\Theta }$. This implies that $G\left( \beta \right)
_{b_{\Theta }}$ contains the identity component $\left( P_{\beta }\right)
_{0}$ of the parabolic subgroup $P_{\beta }=\mathrm{Norm}_{G\left( \beta
\right) }\mathfrak{p}_{\beta }\subset G\left( \beta \right) $. If $\mathfrak{%
g}\left( \beta \right) _{b_{\Theta }}=\mathfrak{p}_{\beta }$ then $G\left(
\beta \right) _{b_{\Theta }}$ is a union of connected components of $%
P_{\beta }$. But $P_{\beta }=Z\left( G(\alpha )\right) (P_{\beta })_{0}$ and 
$Z\left( G(\alpha )\right) \subset M\subset P_{\Theta }$ by Proposition \ref%
{propcentroreal}. Hence, in this case $G\left( \beta \right) _{b_{\Theta
}}=P_{\beta }$ and the orbit $G\left( \beta \right) \cdot b_{\Theta
}=G\left( \beta \right) /P_{\beta }\approx S^{1}$. On the other hand if $%
\mathfrak{p}_{\beta }\neq \mathfrak{g}\left( \beta \right) _{b_{\Theta }}$
then $\mathfrak{g}\left( \beta \right) _{b_{\Theta }}=\mathfrak{g}\left(
\beta \right) $, so that $G\left( \beta \right) _{b_{\Theta }}=G\left( \beta
\right) $ and the orbit $G\left( \beta \right) \cdot b_{\Theta }$ reduces to
a point.

Now, if $\beta \notin \langle \Theta \rangle $ then $\mathfrak{g}_{-\beta }$
has empty intersection with the isotropy subalgebra $\mathfrak{p}_{\Theta }$
which implies that $\mathfrak{g}\left( \beta \right) _{b_{\Theta }}=%
\mathfrak{p}_{\beta }$. Hence $G\left( \beta \right) \cdot b_{\Theta }$ is a
circle $S^{1}$ if $\beta \notin \langle \Theta \rangle $. The last statement
is a consequence of Lemma \ref{lemsuport}.
\end{profe}

We proceed now to look at the orbits of $G(\alpha )$ through $y\in \mathrm{%
att}_{\Theta }(H_{\alpha })$. We recall that $\mathrm{att}_{\Theta
}(H_{\alpha })=\left( K_{H_{\alpha }}\right) _{0}\cdot b_{\Theta }$, where $%
\left( K_{H_{\alpha }}\right) _{0}$ is the identity component of the compact
part of the centralizer $Z_{H_{\alpha }}$ of $H_{\alpha }$. The Lie algebra $%
\mathfrak{k}_{H_{\alpha }}$ of $K_{H_{\alpha }}$ is spanned by $A_{\alpha
}=X_{\alpha }-X_{-\alpha }$.

Let $y=u\cdot b_{\Theta }\in \mathrm{att}_{\Theta }(H_{\alpha })$ with $u\in
\left( K_{H_{\alpha }}\right) _{0}$. Then 
\begin{equation*}
G(\alpha )\cdot y=u\left( u^{-1}G(\alpha )u\right) \cdot b_{\Theta }.
\end{equation*}
The group $u^{-1}G(\alpha )u$ is isomorphic to $G(\alpha )$ and its Lie
algebra $\mathrm{Ad}(u^{-1})\mathfrak{g}(\alpha )$ is isomorphic to $%
\mathfrak{g}(\alpha )$ and hence to $\mathfrak{sl}(2,\mathbb{R})$. Since $%
\mathrm{Ad}\left( u\right) H_{\alpha }=H_{\alpha }$, the decomposition of $%
\mathrm{Ad}(u^{-1})\mathfrak{g}(\alpha )$ into root spaces is given by 
\begin{equation*}
\mathrm{Ad}\left( u^{-1}\right) \mathfrak{g}(\alpha )=\langle H_{\alpha
}\rangle \oplus \mathrm{Ad}\left( u^{-1}\right) \mathfrak{g}_{\alpha }\oplus 
\mathrm{Ad}\left( u^{-1}\right) \mathfrak{g}_{-\alpha }.
\end{equation*}
The subspace $\mathfrak{p}_{u}=\langle H_{\alpha }\rangle \oplus \mathrm{Ad}%
(u^{-1})\mathfrak{g}_{\alpha }$ is a parabolic subalgebra of $\mathrm{Ad}%
(u^{-1})\mathfrak{g}(\alpha )$. Denote by $P_{u}$ the corresponding
parabolic subgroup.

\begin{lema}
\label{isotropia}The subgroup $P_{u}$ is contained in the isotropy subgroup
at $b_{\Theta }$ of the action of $u^{-1}G(\alpha )u$.
\end{lema}

\begin{profe}
We have that $\langle H_{\alpha }\rangle $ is contained in the isotropy
algebra, because $b_{\Theta }$ is a fixed point of $e^{tH_{\alpha }}$. To
see that the same occurs with $\mathrm{Ad}(u^{-1})(\mathfrak{g}_{\alpha })$,
note that if $0\neq X\in \mathrm{Ad}(u^{-1})(\mathfrak{g}_{\alpha })$ then $%
\mathrm{ad}(H_{\alpha })(X)=\alpha (H_{\alpha })X$ because $u$ centralizes $%
H_{\alpha }$. Hence $X$ is an eigenvector of $\mathrm{ad}(H_{\alpha })$
associated to the eigenvalue $\alpha (H_{\alpha })=\langle \alpha ,\alpha
\rangle >0$. Since $H_{\alpha }\in \mathrm{cl}\mathfrak{a}^{+}$ we have that
the eigenspaces of $\mathrm{ad}(H_{\alpha })$ associated to positive
eigenvalues are contained in $\mathfrak{n}^{+}$. Therefore $X\in \mathfrak{n}%
^{+}$ so that $\mathrm{Ad}(u^{-1})(\mathfrak{g}_{\alpha })\subset \mathfrak{n%
}^{+}$, implying that $\mathrm{Ad}(u^{-1})(\mathfrak{g}_{\alpha })$ is
contained in the isotropy subalgebra at $b_{\Theta }$.

Now, the proof follows as in Lemma \ref{lemesse1} by checking that any
connected component of $P_{u}$ is contained in the isotropy subgroup at $%
b_{\Theta }$.
\end{profe}

As a complement to this lemma we have the following homotopy property. From
it the proof of Theorem \ref{teoOrbitsAtt} follows quickly.

\begin{lema}
\label{homotopia}The orbit $\left( u^{-1}G(\alpha )u\right) \cdot b_{\Theta
} $ is a circle $S^{1}$ \ homotopic to the orbit $G(\alpha )\cdot b_{\Theta
} $, by a homotopy that fixes $b_{\Theta }$.
\end{lema}

\begin{profe}
Let $u_{t}\in \left( K_{H_{\alpha }}\right) _{0}$, $t\in \lbrack 0,1]$, be a
continuous curve with $u_{0}=1$ e $u_{1}=u$. Define the continuous map $\psi
:[0,1]\times G(\alpha )\rightarrow \mathbb{F}_{\Theta }$ by 
\begin{equation*}
\psi (t,g)=u_{t}^{-1}gu_{t}\cdot b_{\Theta }.
\end{equation*}
This map has the factorization 
\begin{equation*}
\begin{array}{ccc}
\lbrack 0,1]\times G(\alpha ) & \longrightarrow & \mathbb{F}_{\Theta } \\ 
\downarrow & \nearrow &  \\ 
\left[ 0,1\right] \times \left( G\left( \alpha \right) /P_{\alpha }\right) & 
& 
\end{array}%
\end{equation*}
that defines a continuous map $\phi :[0,1]\times \left( G(\alpha )/P_{\alpha
}\right) \rightarrow \mathbb{F}_{\Theta }$ by $\phi (t,gP_{\alpha })=\psi
(t,g)$. Indeed, $\phi $ is well defined because if $h\in P_{\alpha }$ then $%
u_{t}^{-1}hu_{t}\in P_{u_{t}}$ and, by the previous lemma, $%
(u_{t}^{-1}hu_{t})\cdot b_{\Theta }=b_{\Theta }$. Hence 
\begin{equation*}
\psi (gh)=u_{t}^{-1}gu_{t}(u_{t}^{-1}hu_{t})\cdot b_{\Theta
}=u_{t}^{-1}gu_{t}\cdot b_{\Theta }=\phi \left( g\right)
\end{equation*}
and so $\phi $ is well defined. The function $\phi $ is continuous and if $%
b_{\alpha }$ denotes the origin of $G(\alpha )/P_{\alpha }$ then, again by
the previous lemma, we have that $\phi (t,b_{\alpha })=b_{\Theta }$ for all $%
t\in \lbrack 0,1]$.

Therefore, looking at $G(\alpha )/P_{\alpha }$ as a circle $S^{1}$ with
distinguished point $b_{\alpha }=1\cdot P_{\alpha }$ we see that $\phi $ is
a homotopy between $\phi (0,\cdot )$ whose image is $G(\alpha )\cdot
b_{\Theta }$ and $\phi (1,\cdot )$ whose image is $\left( u^{-1}G(\alpha
)u\right) \cdot b_{\Theta }$. This homotopy fixes $b_{\Theta }$.
\end{profe}

Now we can finish the proof of Theorem \ref{teoOrbitsAtt}: Let $u\in \left(
K_{H_{\alpha }}\right) _{0}$ be such that $y=ub_{\Theta }$. Then $G(\alpha
)\cdot y=u\left( u^{-1}G(\alpha )u\right) \cdot b_{\Theta }$, so that $%
G(\alpha )\cdot y$ and $\left( u^{-1}G(\alpha )u\right) \cdot b_{\Theta }$
are homotopic by a homotopy defined by a curve $t\in \left[ 0,1\right]
\mapsto u_{t}\in \left( K_{H_{\alpha }}\right) _{0}$ with $u_{0}=1$ and $%
u_{1}=u$. By the previous lemma it follows that $G(\alpha )\cdot y$ and $%
G(\alpha )\cdot b_{\Theta }$ are homotopic.

Combining Theorem \ref{teoOrbitsAtt} with Proposition \ref%
{propSemiSimpleGeral} we arrive at once at the following result.

\begin{teorema}
\label{teohomotopyorbits}Let $\alpha $ be a root and $\mathfrak{a}^{+}$ a
Weyl chamber with $H_{\alpha }\in \mathrm{cl}\mathfrak{a}^{+}$. Define the
standard parabolic subgroups $P_{\Theta }$ with respect the positive roots
associated to $\mathfrak{a}^{+}$ and denote by $b_{\Theta }$ the origin of $%
\mathbb{F}_{\Theta }=G/P_{\Theta }$. Take a subset $\Gamma \subset G$ with $%
G\left( \alpha \right) \subset \Gamma $. Assume that

\begin{enumerate}
\item for every minimal flag manifold $\mathbb{F}_{\Theta }$ the orbit $%
G\left( \alpha \right) \cdot b_{\Theta }$ (which is a closed curve) is not
homotopic to a point and

\item the subgroup generated by $\Gamma $ is Zariski dense.
\end{enumerate}

Then $\Gamma $ generates $G$ as semigroup.
\end{teorema}

As will be seen below condition (1) of this theorem holds in three cases
namely when $\mathfrak{g}=\mathfrak{sl}\left( n,\mathbb{R}\right) $, $\alpha 
$ is a long root of $\mathfrak{sp}\left( n,\mathbb{R}\right) $ and $\alpha $
is a short root of the $G_{2}$ digram.

\subsection{The fundamental group of minimal flag manifolds}

We look here at the fundamental groups of the minimal flag manifolds and the
homotopy classes of the orbits $G\left( \alpha \right) \cdot b_{\Theta }$
appearing in Theorem \ref{teohomotopyorbits}.

Fix a simple system of roots $\Sigma $. For a root $\alpha $ we choose $%
X_{\alpha }\in \mathfrak{g}_{\alpha }$, $X_{-\alpha }=-\theta X_{\alpha }\in 
\mathfrak{g}_{-\alpha }$ such that $\langle X_{\alpha },X_{-\alpha }\rangle
=1$, and write $A_{\alpha }=X_{\alpha }-X_{-\alpha }\in \mathfrak{k}_{\alpha
}$ and $S_{\alpha }=X_{\alpha }+X_{-\alpha }$. If $\mathbb{F}_{\Theta }$, $%
\Theta \subset \Sigma $, is a flag manifold we write $b_{\Theta }$ for its
origin defined by $\Sigma $. Then for a root $\alpha $ the orbit $G\left(
\alpha \right) \cdot b_{\Theta }$ is the closed curve $e^{tA_{\alpha }}\cdot
b_{\Theta }$. We denote by $c_{\alpha }^{\Theta }$ (or simply $c_{\alpha }$)
the homotopy class of this curve in the fundamental group of $\mathbb{F}%
_{\Theta }$.

The fundamental groups of real flag manifolds were described in Wiggerman 
\cite{wig} by generators and relations. We recall the result of \cite{wig}
the case of a split real form $\mathfrak{g}$.

\begin{teorema}
\label{teofundwig}The fundamental group $\pi _{1}\left( \mathbb{F}_{\Theta
}\right) $ of $\mathbb{F}_{\Theta }$ is generated by $c_{\alpha }$ with $%
\alpha \in \Sigma $, subjected to the relations

\begin{enumerate}
\item $c_{\alpha }=1$ if $\alpha \in \Theta $ and

\item $c_{\alpha }c_{\beta }c_{\alpha }^{-1}c_{\beta }^{-\varepsilon \left(
\alpha ,\beta \right) }=1$ where $\varepsilon \left( \alpha ,\beta \right)
=\left( -1\right) ^{\langle \alpha ^{\vee },\beta \rangle }$ and $\langle
\alpha ^{\vee },\beta \rangle =\frac{2\langle \alpha ^{\vee },\beta \rangle 
}{\langle \alpha ,\alpha \rangle }$ is the Killing number.
\end{enumerate}
\end{teorema}

The first relation says that $\pi _{1}\left( \mathbb{F}_{\Theta }\right) $
is in fact generated by $c_{\alpha }$ with $\alpha \in \Sigma \setminus
\Theta $. However it is convenient to include in the statement the
generators $c_{\alpha }$, $\alpha \in \Theta $, because they enter in the
second set of relations. Namely, $c_{\beta }c_{\beta }^{-\varepsilon \left(
\alpha ,\beta \right) }=c_{\alpha }c_{\beta }c_{\alpha }^{-1}c_{\beta
}^{-\varepsilon \left( \alpha ,\beta \right) }=1$ if $\alpha \in \Theta $
and $\beta \in \Sigma \setminus \Theta $. Notice that this relation implies
that $c_{\beta }^{2}=1$, $\beta \in \Sigma \setminus \Theta $, if there
exists a root $\alpha \in \Theta $ such that the Killing number $\langle
\alpha ^{\vee },\beta \rangle $ is odd.

From these generators and relations it is easy to get the fundamental groups
of the minimal flag manifolds. Given a root $\beta \in \Sigma $ we write $%
\Theta _{\beta }=\Sigma \setminus \{\beta \}$ and take the corresponding
minimal flag manifold $\mathbb{F}_{\Theta _{\beta }}$. These exaust the
minimal flag manifolds, except for the diagram $A_{1}$ when the only flag
manifold is given by $\Theta =\emptyset $.

\begin{proposicao}
\label{propze2}The fundamental group of a minimal flag manifold $\mathbb{F}%
_{\Theta _{\beta }}$ is $\pi _{1}\left( \mathbb{F}_{\Theta _{\beta }}\right)
=\mathbb{Z}_{2}$ except in for $A_{1}$ or when $\beta $ is the long root in
the $C_{l}$ diagram. In the exceptions $\pi _{1}\left( \mathbb{F}_{\Theta
_{\beta }}\right) =\mathbb{Z}$.
\end{proposicao}

\begin{profe}
By Theorem \ref{teofundwig} the fundamental group of $\mathbb{F}_{\Theta
_{\beta }}$ is cyclic and generated by $c_{\beta }$. If the diagram is not $%
A_{1}$ then $\beta $ is linked to another root $\alpha $. If the link is a
simple edge then $\langle \alpha ^{\vee },\beta \rangle =-1$ and hence $%
c_{\beta }^{2}=1$ so that $\pi _{1}\left( \mathbb{F}_{\Theta _{\beta
}}\right) =\mathbb{Z}_{2}$. The Killing number $\langle \alpha ^{\vee
},\beta \rangle $ is also odd in the $G_{2}$ diagram ($-1$ or $-3$) or in
case $\alpha $ is a long root and $\beta $ a short root. A glance at the
Dynkin diagrams shows that every root $\beta $ has such link if $\beta $ is
not the long root in the $C_{l}$ diagram. If $\beta $ is this long root then
there are no relations involving $c_{\beta }$ and hence $\pi _{1}\left( 
\mathbb{F}_{\Theta _{\beta }}\right) $ is cyclic infinity.

Finally, the flag manifold of $A_{1}$ is the circle $S^{1}$ and hence has $%
\mathbb{Z}$ as fundamental group.
\end{profe}

We proceed now to look at the homotopy condition of Theorems \ref%
{teoOrbitsAtt} and \ref{teohomotopyorbits}. Here we change slightly the
point of view. In those theorems we started with a root $\alpha $ and
choosed a Weyl chamber containing $H_{\alpha }$ in its closure. Here we fix
a Weyl chamber $\mathfrak{a}^{+}$ (and hence an origin $b_{\Theta }$ in a
flag manifold $\mathbb{F}_{\Theta }$) and take a root $\mu $ such that $%
H_{\mu }\in \mathrm{cl}\mathfrak{a}^{+}$. Since the Weyl group $\mathcal{W}$
acts transitively on the set of Weyl chambers, there is no loss of
generality in fixing $\mathfrak{a}^{+}$ in advance. In fact, by $wH_{\alpha
}=H_{w\alpha }$ it follows that for any root $\alpha $ there exists a root $%
\mu =w\alpha $ such that $H_{\mu }\in \mathrm{cl}\mathfrak{a}^{+}$, that is, 
$H_{\alpha }\in \mathrm{cl}\left( w^{-1}\mathfrak{a}^{+}\right) $ . If $%
\overline{w}\in M^{\ast }$ is a representative of $w$ then $G\left( \alpha
\right) =\overline{w}^{-1}G\left( \mu \right) \overline{w}$ and $G\left(
\alpha \right) \cdot \overline{w}^{-1}b_{\Theta }=\overline{w}^{-1}G\left(
\mu \right) \overline{w}\cdot \overline{w}^{-1}b_{\Theta }=\overline{w}%
^{-1}\left( G\left( \mu \right) \cdot b_{\Theta }\right) $. The point $%
\overline{w}^{-1}b_{\Theta }$ is the origin corresponding to $w^{-1}%
\mathfrak{a}^{+}$ and a curve $t\in \left[ 0,1\right] \mapsto g_{t}\in G$
with $g_{0}=1$ and $g_{1}=\overline{w}$ realizes a homotopy between the two
orbits $G\left( \alpha \right) \cdot \overline{w}^{-1}b_{\Theta }$ and $%
G\left( \mu \right) \cdot b_{\Theta }$. Therefore we can restrict our
analysis to roots $\mu $ with $H_{\mu }\in \mathrm{cl}\mathfrak{a}^{+}$ and
the Weyl chamber $\mathfrak{a}^{+}$ previously fixed.

The action of $\mathcal{W}$ on the set of roots is either transitive (for
the simply laced diagrams $A_{l}$, $D_{l}$, $E_{6}$, $E_{7}$ and $E_{8}$) or
has two orbits the long and short roots (for the other diagrams $B_{l}$, $%
C_{l}$, $F_{4}$ and $G_{2}$). This implies that for the simply laced
diagrams there is just one root (the highest positive root) $\mu $ with $%
H_{\mu }\in \mathrm{cl}\mathfrak{a}^{+}$ while in the other cases there is
the highest root and a short root as well in the closure of the Weyl chamber.

To look at the classes $c_{\mu }^{\Theta }$ with $H_{\mu }\in \mathrm{cl}%
\mathfrak{a}^{+}$ in the fundamental groups of the minimal flag manifolds we
realize them as orbits of projective representations.

\subsection{Projective and spherical orbits}

To look at the fundamental group of the flag manifolds we shall exploit
their realizations as projective orbits of representations. Given the simple
system of roots $\Sigma =\{\alpha _{1},\ldots ,\alpha _{l}\}$ let $\{\omega
_{1},\ldots ,\omega _{l}\}$ be the set of basic weights defined by 
\begin{equation*}
\langle \alpha _{i}^{\vee },\omega _{j}\rangle =\frac{2\langle \alpha
_{i},\omega _{j}\rangle }{\langle \alpha _{i},\alpha _{j}\rangle }=\delta
_{ij}.
\end{equation*}
Any $\omega =p_{1}\omega _{1}+\cdots +p_{l}\omega _{l}$, $p_{i}\in \mathbb{Z}
$, $p_{i}\geq 0$, is the highest weight of a representation $\rho _{\omega }$
of $\mathfrak{g}$ in a vector space $V_{\omega }$. Write $G=\langle \exp
\rho _{\omega }\left( \mathfrak{g}\right) \rangle $ for the group with Lie
algebra $\mathfrak{g}$ that integrates the representation.

Let $V\left( \omega \right) \subset V_{\omega }$ be the one-dimensional
weight space of $\omega $. For $v\in V\left( \omega \right) $, $v\neq 0$,
put $V^{+}\left( \omega \right) =\mathbb{R}^{+}v$ for the ray spanned by $v$.

Consider the projective orbit $G\cdot V\left( \omega \right) =\{g\cdot
V\left( \omega \right) :g\in G\}$ and the spherical orbit $G\cdot
V^{+}\left( \omega \right) =\{g\cdot V^{+}\left( \omega \right) :g\in G\}$.
The isotropy subalgebra at both $V\left( \omega \right) $ and $V^{+}\left(
\omega \right) $ is the parabolic subalgebra $\mathfrak{p}_{\Theta _{\omega
}}$ where 
\begin{equation*}
\Theta _{\omega }=\{\alpha _{i}\in \Sigma :p_{i}=0\}=\{\alpha \in \Sigma
:\alpha \left( \omega \right) =0\}.
\end{equation*}
Hence the identity component of the parabolic subgroup $P_{\Theta _{\omega
}} $ is contained in the isotropy subgroups of $V\left( \omega \right) $ and 
$V^{+}\left( \omega \right) $. It turns out that the isotropy subgroup at $%
V\left( \omega \right) $ is the whole $P_{\Theta _{\omega }}$, since $m\cdot
v=\pm v$ if $m\in M$. Therefore $G\cdot V\left( \omega \right) \approx 
\mathbb{F}_{\Theta _{\omega }}$ with the origin $b_{\Theta _{\omega }}\in 
\mathbb{F}_{\Theta _{\omega }}$ being identified to $V\left( \omega \right)
\in G\cdot V\left( \omega \right) $. Also the the map $G\cdot V^{+}\left(
\omega \right) \rightarrow G\cdot V\left( \omega \right) $, $g\cdot
V^{+}\left( \omega \right) \mapsto g\cdot V\left( \omega \right) $ is a
covering because the isotropy subgroup at $V^{+}\left( \omega \right) $ is
an open subgroup of $P_{\Theta _{\omega }}$. This covering has one or two
leaves, depending if $-V^{+}\left( \omega \right) $ belongs or not to the
orbit $G\cdot V^{+}\left( \omega \right) $.

To look closely at the covering $G\cdot V^{+}\left( \omega \right)
\rightarrow G\cdot V\left( \omega \right) $ we take a root $\alpha $. The
orbit $G\left( \alpha \right) \cdot b_{\Theta }$ is the curve $\left( \exp
tA_{\alpha }\right) \cdot b_{\Theta }$, $t\in \left[ 0,\pi \right] $, where $%
A_{\alpha }=X_{\alpha }-X_{-\alpha }$. This means that $m_{\alpha }v=\pm v$
if $v\in V\left( \omega \right) $, where $m_{\alpha }=\exp \pi \rho _{\omega
}\left( A_{\alpha }\right) $. The following lemma gives the sign in this
equality.

\begin{lema}
\label{lemmpesomax}$m_{\alpha }v=\left( -1\right) ^{\omega \left( H_{\alpha
}^{\vee }\right) }v$.
\end{lema}

\begin{profe}
The Lie algebra $\rho _{\omega }\left( \mathfrak{g}\right) $ as well as the
Lie group $G=\langle \exp \rho _{\omega }\left( \mathfrak{g}\right) \rangle $
can be complexified. In the complexification we have $m_{\alpha }=\exp \pi
\rho _{\omega }\left( A_{\alpha }\right) =\exp i\pi \rho _{\omega }\left(
H_{\alpha }^{\vee }\right) $ (cf. Proposition \ref{propcentroreal}). Hence $%
m_{\alpha }v=e^{i\pi \omega \left( H_{\alpha }^{\vee }\right) }v=\left(
-1\right) ^{\omega \left( H_{\alpha }^{\vee }\right) }v$ if $v\in V\left(
\omega \right) $.
\end{profe}

Now let $\omega =\omega _{j}$ be a basic weight, so that $\Theta _{\omega
_{j}}=\Sigma \setminus \{\alpha _{j}\}$, and $\mathbb{F}_{\Theta _{\omega
_{j}}}$ is a minimal flag manifold. By the very definition of the basic
weights we have $\omega _{j}\left( H_{\alpha _{j}}^{\vee }\right) =1$, so
that $m_{\alpha _{j}}v=-v$ if $v\in V\left( \omega _{j}\right) $. This
implies that the spherical orbit $G\cdot V^{+}\left( \omega _{j}\right) $ is
a double covering of $\mathbb{F}_{\Theta _{\omega _{j}}}$.

On the other hand we have, by Proposition \ref{propze2}, that except in two
cases the fundamental group of a minimal flag manifold is $\mathbb{Z}_{2}$.
Combining these facts we get the universal covering space of the minimal
flag manifolds.

\begin{proposicao}
\label{propunivercover}Let $\omega _{j}$ be a basic weight in a diagram
different from $A_{1}$ and such that $\alpha _{j}$ is not the long root of
the $C_{l}$ diagram. Let $\mathbb{F}_{\Theta _{\omega _{j}}}$ be the
corresponding minimal flag manifold. Then the spherical orbit $G\cdot
V^{+}\left( \omega _{j}\right) $ is the simply connected universal cover of $%
\mathbb{F}_{\Theta _{\omega _{j}}}$.
\end{proposicao}

\begin{profe}
In fact, $\mathbb{F}_{\Theta _{\omega _{j}}}$ has only two coverings since
its fundamental group is $\mathbb{Z}_{2}$. Since $G\cdot V^{+}\left( \omega
_{j}\right) \rightarrow \mathbb{F}_{\Theta _{\omega _{j}}}$ is a double
covering it follows that $G\cdot V^{+}\left( \omega _{j}\right) $ is indeed
the simply connected cover.
\end{profe}

As a consequence we have that a closed curve in $\mathbb{F}_{\Theta _{\omega
_{j}}}$ is null homotopic if and only if its lifting to $G\cdot V^{+}\left(
\omega _{j}\right) $ is a closed curve. This yields the following criterion
to decide if an orbit $G\left( \alpha \right) \cdot b_{\Theta }$ is null
homotopic.

\begin{proposicao}
Let $\mathbb{F}_{\Theta _{\omega _{j}}}$ be a minimal flag manifold with $%
\omega _{j}$ as in Proposition \textrm{{\ref{propunivercover}}}. If $\alpha $
is a root then $G\left( \alpha \right) \cdot b_{\Theta _{\omega _{j}}}$ is
null homotopic if and only if the Killing number $\omega _{j}\left(H_{\alpha
}^{\vee }\right) $ is even.
\end{proposicao}

\begin{profe}
The orbit $G\left( \alpha \right) \cdot b_{\Theta _{\omega _{j}}}$ is the
closed curve $\left( \exp tA_{\alpha }\right) \cdot b_{\Theta _{\omega
_{j}}} $, $t\in \left[ 0,\pi \right] $. In terms of the projective orbit of
the representation this curve is given by $\exp t\rho _{\omega _{j}}\left(
A_{\alpha }\right) \cdot V\left( \omega _{j}\right) $, $t\in \left[ 0,\pi %
\right] $. The lifting of this curve to $G\cdot V^{+}\left( \omega
_{j}\right) $ starting at $V^{+}\left( \omega _{j}\right) $ is $\exp t\rho
_{\omega _{j}}\left( A_{\alpha }\right) \cdot V^{+}\left( \omega _{j}\right) 
$, $t\in \left[ 0,\pi \right] $. This lifting is a closed curve if and only
if $m_{\alpha }\cdot v=\exp \pi \rho _{\omega _{j}}\left( A_{\alpha }\right)
\cdot v=v$ for $v\in V\left( \omega _{j}\right) $. By Lemma \ref{lemmpesomax}
we have $m_{\alpha }v=\left( -1\right) ^{\omega _{j}\left( H_{\alpha }^{\vee
}\right) }v$, so that $G\left( \alpha \right) \cdot b_{\Theta _{\omega
_{j}}} $ is null homotopic if and only if $\omega _{j}\left( H_{\alpha
}^{\vee }\right) $ is even.
\end{profe}

If $\alpha =n_{1}\alpha _{1}+\cdots +n_{l}\alpha _{l}$ where $\Sigma
=\{\alpha _{1},\ldots ,\alpha _{l}\}$ is the simple system of roots then we
can compute $\omega _{j}\left( H_{\alpha }^{\vee }\right) $ explicitly from
the coefficients $n_{j}$. In fact, since $\langle \alpha _{i},\omega
_{j}\rangle =0$ if $i\neq j$ and $2\langle \alpha _{j},\omega _{j}\rangle
=\langle \alpha _{j},\alpha _{j}\rangle $, we get 
\begin{equation}
\omega _{j}\left( H_{\alpha }^{\vee }\right) =\frac{2\langle \alpha ,\omega
_{j}\rangle }{\langle \alpha ,\alpha \rangle }=n_{j}\frac{2\langle \alpha
_{j},\omega _{j}\rangle }{\langle \alpha ,\alpha \rangle }=n_{j}\frac{%
\langle \alpha _{j},\alpha _{j}\rangle }{\langle \alpha ,\alpha \rangle }.
\label{fornumkillingpesoraiz}
\end{equation}

By this formula we see that $G\left( \alpha \right) \cdot b_{\Theta _{\omega
_{j}}}$ is null homotopic if $n_{j}$ is even and the roots have the same
length. We observe the following possibilities when the roots $\alpha _{j}$
and $\alpha $ have different length.

\begin{enumerate}
\item $\alpha _{j}$ is long and $\alpha $ is short with $\langle \alpha
_{j},\alpha _{j}\rangle =2\langle \alpha ,\alpha \rangle $. Then $G\left(
\alpha \right) \cdot b_{\Theta _{\omega _{j}}}$ is null homotopic regardless
the coefficient $n_{j}$.

\item $\alpha $ is long and $\alpha _{j}$ is short with $\langle \alpha
,\alpha \rangle =2\langle \alpha _{j},\alpha _{j}\rangle $. Then $G\left(
\alpha \right) \cdot b_{\Theta _{\omega _{j}}}$ is not null homotopic if $%
n_{j}\neq 4p$.
\end{enumerate}

\subsection{Roots in the closure of the Weyl chamber}

\label{sectionclosure}

Now we can look at the homotopy class of the orbits $G\left( \alpha \right)
\cdot b_{\Theta }$ in a minimal flag manifold, for a root $\alpha $ with $%
H_{\alpha }\in \mathrm{cl}\mathfrak{a}^{+}$. Concerning the highest roots we
write their coefficients above the simple roots in the Dynkin diagrams:

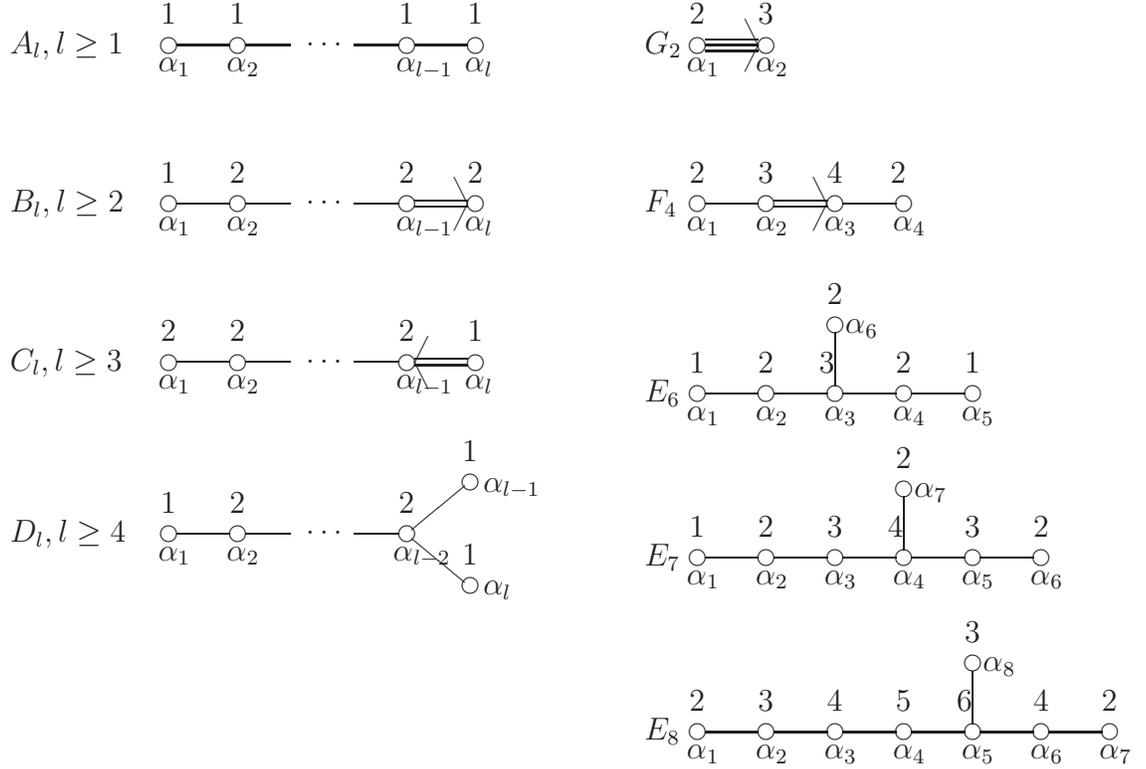
\begin{picture}(400,330)(-20,0)

\put(-20,260){
\begin{picture}(140,60)(-40,0)

\put(-40,27){$A_l,l \geq 1$}
\put(20,30){\circle{6}}
\put(17,38){$1$}
\put(23,30){\line(1,0){20}}
\put(46,30){\circle{6}}
\put(43,38){$1$}

\put(49,30){\line(1,0){17}}
\put(72,30){\ldots}
\put(90,30){\line(1,0){17}}
\put(110,30){\circle{6}}
\put(107,38){$1$}

\put(113,30){\line(1,0){20}}
\put(136,30){\circle{6}}
\put(133,38){$1$}

\put(16,20){${\alpha}_1$}
\put(42,20){${\alpha}_2$}
\put(106,20){${\alpha}_{l-1}$}
\put(132,20){${\alpha}_l$}
\end{picture}
}

\put(-20,200){
\begin{picture}(100,60)(-40,0)

\put(-40,27){$B_l,l \geq 2$}
\put(20,30){\circle{6}}
\put(17,38){$1$}

\put(23,30){\line(1,0){20}}
\put(46,30){\circle{6}}
\put(43,38){$2$}

\put(49,30){\line(1,0){17}}
\put(72,30){\ldots}
\put(90,30){\line(1,0){17}}
\put(110,30){\circle{6}}
\put(107,38){$2$}

\put(113,31.2){\line(1,0){20}}
\put(113,28.8){\line(1,0){20}}
\put(136,30){\circle{6}}
\put(133,38){$2$}

\put(133,30){\line(-1,2){5}}
\put(133,30){\line(-1,-2){5}}
\put(16,20){${\alpha}_1$}
\put(42,20){${\alpha}_2$}
\put(106,20){${\alpha}_{l-1}$}
\put(132,20){${\alpha}_l$}

\end{picture}
}
\put(-20,140){
\begin{picture}(100,60)(-40,0)

\put(-40,27){$C_l,l \geq 3$}
\put(20,30){\circle{6}}
\put(17,38){$2$}

\put(23,30){\line(1,0){20}}
\put(46,30){\circle{6}}
\put(43,38){$2$}

\put(49,30){\line(1,0){17}}
\put(72,30){\ldots}
\put(90,30){\line(1,0){17}}
\put(110,30){\circle{6}}
\put(107,38){$2$}

\put(113,30){\line(1,2){5}}
\put(113,30){\line(1,-2){5}}
\put(113,31.2){\line(1,0){20}}
\put(113,28.8){\line(1,0){20}}
\put(136,30){\circle{6}}
\put(133,38){$1$}

\put(16,20){${\alpha}_1$}
\put(42,20){${\alpha}_2$}
\put(106,20){${\alpha}_{l-1}$}
\put(132,20){${\alpha}_l$}
\end{picture}
}
\put(-20,70){
\begin{picture}(100,70)(-40,0)

\put(-40,32){$D_l,l \geq 4$}
\put(20,35){\circle{6}}
\put(17,43){$1$}

\put(16,25){$\alpha_1$}
\put(23,35){\line(1,0){20}}
\put(46,35){\circle{6}}
\put(43,43){$2$}

\put(42,25){$\alpha_2$}
\put(49,35){\line(1,0){17}}
\put(72,35){\ldots}
\put(90,35){\line(1,0){17}}
\put(110,35){\circle{6}}
\put(107,43){$2$}

\put(105,25){$\alpha_{l-2}$}
\put(111.9,36.6){\line(6,5){20}}
\put(111.9,33.4){\line(6,-5){20}}
\put(134,54.6){\circle{6}}
\put(131,62.6){$1$}

\put(139,51){$\alpha_{l-1}$}
\put(134,15.6){\circle{6}}
\put(131,23.6){$1$}

\put(139,12){$\alpha_l$}

\end{picture}
}
\put(220,260){
\begin{picture}(100,60)(0,0)

\put(0,27){$G_2$}
\put(20,30){\circle{6}}
\put(17,38){$2$}

\put(23,30){\line(1,0){20}}
\put(23,32.2){\line(1,0){20}}
\put(23,27.8){\line(1,0){20}}
\put(46,30){\circle{6}}
\put(43,38){$3$}

\put(43,30){\line(-1,2){5}}
\put(43,30){\line(-1,-2){5}}
\put(16,20){${\alpha}_1$}
\put(42,20){${\alpha}_2$}
\end{picture}
}
\put(220,200){
\begin{picture}(100,60)(0,0)

\put(0,27){$F_4$}
\put(20,30){\circle{6}}
\put(17,38){$2$}

\put(16,20){${\alpha}_1$}
\put(23,30){\line(1,0){20}}
\put(46,30){\circle{6}}
\put(43,38){$3$}

\put(42,20){${\alpha}_2$}
\put(49,31.2){\line(1,0){20}}
\put(49,28.8){\line(1,0){20}}
\put(72,30){\circle{6}}
\put(69,38){$4$}

\put(68,20){${\alpha}_3$}
\put(69,30){\line(-1,2){5}}
\put(69,30){\line(-1,-2){5}}
\put(75,30){\line(1,0){20}}
\put(98,30){\circle{6}}
\put(93,38){$2$}

\put(94,20){${\alpha}_4$}
\end{picture}
}
\put(220,128){
\begin{picture}(130,70)(0,0)

\put(0,27){$E_6$}
\put(20,30){\circle{6}}
\put(17,38){$1$}

\put(23,30){\line(1,0){20}}
\put(46,30){\circle{6}}
\put(43,38){$2$}

\put(49,30){\line(1,0){20}}
\put(72,30){\circle{6}}
\put(66,38){$3$}

\put(75,30){\line(1,0){20}}
\put(98,30){\circle{6}}
\put(95,38){$2$}

\put(101,30){\line(1,0){20}}
\put(124,30){\circle{6}}
\put(121,38){$1$}

\put(72,33){\line(0,1){20}}
\put(72,56){\circle{6}}
\put(69,64){$2$}

\put(16,20){${\alpha}_1$}
\put(42,20){${\alpha}_2$}
\put(68,20){${\alpha}_3$}
\put(94,20){${\alpha}_4$}
\put(120,20){${\alpha}_5$}
\put(76,53){${\alpha}_6$}

\end{picture}
}
\put(220,66){
\begin{picture}(200,70)(0,0)

\put(0,27){$E_7$}
\put(20,30){\circle{6}}
\put(17,38){$1$}

\put(23,30){\line(1,0){20}}
\put(46,30){\circle{6}}
\put(43,38){$2$}

\put(49,30){\line(1,0){20}}
\put(72,30){\circle{6}}
\put(69,38){$3$}

\put(75,30){\line(1,0){20}}
\put(98,30){\circle{6}}
\put(92,38){$4$}

\put(101,30){\line(1,0){20}}
\put(124,30){\circle{6}}
\put(121,38){$3$}

\put(127,30){\line(1,0){20}}
\put(150,30){\circle{6}}
\put(147,38){$2$}

\put(98,33){\line(0,1){20}}
\put(98,56){\circle{6}}
\put(95,64){$2$}

\put(16,20){${\alpha}_1$}
\put(42,20){${\alpha}_2$}
\put(68,20){${\alpha}_3$}
\put(94,20){${\alpha}_4$}
\put(120,20){${\alpha}_5$}
\put(146,20){${\alpha}_6$}
\put(102,53){${\alpha}_7$}

\end{picture}
}
\put(220,0){
\begin{picture}(200,70)(0,0)

\put(0,27){$E_8$}
\put(20,30){\circle{6}}
\put(17,38){$2$}

\put(23,30){\line(1,0){20}}
\put(46,30){\circle{6}}
\put(43,38){$3$}

\put(49,30){\line(1,0){20}}
\put(72,30){\circle{6}}
\put(69,38){$4$}

\put(75,30){\line(1,0){20}}
\put(98,30){\circle{6}}
\put(95,38){$5$}

\put(101,30){\line(1,0){20}}
\put(124,30){\circle{6}}
\put(118,38){$6$}

\put(127,30){\line(1,0){20}}
\put(150,30){\circle{6}}
\put(147,38){$4$}

\put(153,30){\line(1,0){20}}
\put(176,30){\circle{6}}
\put(173,38){$2$}

\put(124,33){\line(0,1){20}}
\put(124,56){\circle{6}}
\put(121,64){$3$}

\put(16,20){${\alpha}_1$}
\put(42,20){${\alpha}_2$}
\put(68,20){${\alpha}_3$}
\put(94,20){${\alpha}_4$}
\put(120,20){${\alpha}_5$}
\put(146,20){${\alpha}_6$}
\put(172,20){${\alpha}_7$}
\put(128,53){${\alpha}_8$}

\end{picture}
}

\end{picture}%
%
%

Consider first the highest root $\mu $. By inspecting the coefficients of $%
\mu $ we can apply formula (\ref{fornumkillingpesoraiz}) to get the
following to cases where $G\left( \mu \right) \cdot b_{\Theta }$ is not null
homotopic in whatsoever minimal flag manifold:

\begin{enumerate}
\item $A_{l}$, because the coefficients are $n_{j}=1$ and the roots have the
same length.

\item $C_{l}$, in this case the highest root $\mu $ has coefficients $%
n_{j}=2 $ only with respect to the short roots $\alpha _{j}$, $j<l$. So that 
$\omega _{j}\left( H_{\mu }^{\vee }\right) =1$ for every $j$.
\end{enumerate}

In the diagrams $D_{l}, E_{6}, E_{7}$ and $E_{8} $, the roots have the same
length and the coefficients $n_{j}$ are even for several $j $. This means
that for these diagrams $G\left( \mu \right) \cdot b_{\Theta }$ is null
homotopic on several minimal flag manifolds. We have also that in the
diagrams $B_{l}, G_{2}$ and $F_{4}$, the Killing number $\omega _{j}\left(
H_{\mu }^{\vee }\right) $ is even if $\alpha _{j}$ is a short simple root
for several $j $. Hence for these diagrams $G\left( \mu \right) \cdot
b_{\Theta }$ is null homotopic on several minimal flag manifolds.

Now let $\overline{\mu }$ be the short root with $H_{\overline{\mu }}\in 
\mathrm{cl}\mathfrak{a}^{+}$ in the diagrams with multiple edges. For $B_{l}$%
, $C_{l}$ and $F_{4}$, the Killing number $\omega _{j}\left( H_{\overline{%
\mu }}^{\vee }\right) $ is even if $\alpha _{j}$ is a long simple root.
Hence for these diagrams $G\left( \overline{\mu }\right) \cdot b_{\Theta }$
is null homotopic on several minimal flag manifolds.

On the other hand, in $G_{2}$ we have $\overline{\mu }=\alpha _{1}+3\alpha
_{2}$. Hence by formula (\ref{fornumkillingpesoraiz}), $\omega _{j}\left( H_{%
\overline{\mu }}^{\vee }\right) $ is odd for $j=1,2$. This means that for $%
G_{2}$ the orbits $G\left( \overline{\mu }\right) \cdot b_{\Theta }$ in the
minimal flag manifolds are not null homotopic.

Having these facts about the fundamental groups Theorem \ref%
{teohomotopyorbits} yields immediately following result.

\begin{teorema}
\label{teoreal}Let $\Gamma \subset G$ be a subset such that $G\left( \alpha
\right) \subset \Gamma $ and the group generated by $\Gamma $ is Zariski
dense. Then $\Gamma $ generates $G$ in the following cases:

\begin{enumerate}
\item $\mathfrak{g}=\mathfrak{sl}\left( l+1,\mathbb{R}\right) $.

\item $\mathfrak{g}=\mathfrak{sp}\left( l,\mathbb{R}\right) $ and $\alpha $
is a long root.

\item $\mathfrak{g}$ is the split real form associated to $G_{2}$ and $%
\alpha $ is a short root.
\end{enumerate}
\end{teorema}

\subsection{\protect\medskip Example: Short roots in $\mathfrak{sp}\left( n,%
\mathbb{R}\right) $}

\label{secexampl}

We present an\textbf{\ }example showing that the result of Theorem \ref%
{teoreal} does not hold if $\alpha $ is a short root of the symplectic Lie
algebra $\mathfrak{sp}\left( l,\mathbb{R}\right) $. Let $Q$ be the quadratic
form in $\mathbb{R}^{2l}$ whose matrix with respect to the standard basis $%
\{e_{1},\ldots ,e_{l},f_{1},\ldots ,f_{l}\}$ is 
\begin{equation*}
\lbrack Q]=\left( 
\begin{array}{ll}
0 & \mathrm{id}_{l\times l} \\ 
\mathrm{id}_{l\times l} & 0%
\end{array}
\right).
\end{equation*}
Define the subset of the projective space $C\subset \mathbb{P}^{2l-1}$ by 
\begin{equation*}
C=\{[v]\in \mathbb{P}^{2l-1}:Q\left( v\right) \geq 0\}
\end{equation*}
where $[v]=\mathrm{span}\{v\}$, $0\neq v\in \mathbb{R}^{2l}$. It is easily
seen that $C$ is compact and has interior $\mathrm{int}C=\{[v]\in \mathbb{P}
^{2l-1}:Q\left( v\right) >0\}$.

Now define $S\subset \mathrm{Sp}{\left( l,\mathbb{R}\right) }$ to be the
compression semigroup 
\begin{equation*}
S=\{g\in \mathrm{Sp}\left( l,\mathbb{R}\right) {:gC\subset C\}.}
\end{equation*}
We claim that $\mathrm{int}S\neq \emptyset $. In fact, $X=[Q]$ is itself a
symplectic matrix and for any $0\neq v\in \mathbb{R}^{2l}$ we have 
\begin{equation*}
\frac{d}{dt}Q\left( e^{tX}v\right) =\left( e^{tX}v\right) ^{T}\left(
X^{T}[Q]+[Q]X\right) e^{tX}v=2\left( e^{tX}v\right) ^{T}e^{tX}v>0,
\end{equation*}
so that $t\mapsto Q\left( e^{tX}v\right) $ is strictly increasing. This
implies that for $t>0$, $e^{tX}C\subset \mathrm{int}C\subset C$ hence $%
e^{tX}\in S$ if $t\geq 0$. Actually, $e^{tX}\in \mathrm{int}S$ if $t>0$ by
continuity with respect to the compact-open topology.

Now, take a symplectic matrix of the form 
\begin{equation*}
Y=\left( 
\begin{array}{ll}
A & 0 \\ 
0 & -A^{T}%
\end{array}
\right) .
\end{equation*}
Then $Y^{T}[Q]+[Q]Y=0$, so that $e^{tY}$, $t\in \mathbb{R}$, is an isometry
of $Q$, that is, $Q\left( e^{tY}v\right) =Q\left( v\right) $ for $v\in 
\mathbb{R}^{2l-1}$. Hence $e^{tY}\in S$, $t\in \mathbb{R}$, which in turn
implies that the group 
\begin{equation*}
G=\{\left( 
\begin{array}{ll}
g & 0 \\ 
0 & \left( g^{-1}\right) ^{T}%
\end{array}
\right) \in \mathrm{Sp}\left( l,\mathbb{R}\right) :g\in \mathrm{Gl}
^{+}\left( l,\mathbb{R}\right) \}
\end{equation*}
is contained in $S$. But for any short root $\alpha =\lambda _{i}-\lambda
_{j}$, $i\neq j$, we have $G\left( \alpha \right) \subset G$, concluding our
example.

\section{Representations of $\mathfrak{sl}\left( 2,\mathbb{C}\right) $}

In this section we take an irreducible representation of $\mathfrak{sl}%
\left( 2,\mathbb{C}\right) $ on $\mathbb{C}^{n+1}$, that is, a homomorphism $%
\rho _{n}:\mathfrak{sl}\left( 2,\mathbb{C}\right) \rightarrow \mathfrak{sl}%
\left( n+1,\mathbb{C}\right) $ and look at the subgroup $G_{1}=\langle \exp
\rho _{n}\left( \mathfrak{sl}\left( 2,\mathbb{C}\right) \right) \rangle
\subset \mathrm{Sl}\left( n+1,\mathbb{C}\right) $. As an application of the
previous results will show the following result.

\begin{teorema}
\label{teoParaRepreSl2}Let $\rho _{n}:\mathfrak{sl}\left( 2,\mathbb{C}%
\right) \rightarrow \mathfrak{sl}\left( n+1,\mathbb{C}\right) $ be the
irreducible $\left( n+1\right) $-dimensional representation of $\mathfrak{sl}%
\left( 2,\mathbb{C}\right) $ and suppose that $\Gamma \subset \mathrm{Sl}%
\left( n+1,\mathbb{C}\right) $ is a subset containing $G_{1}=\langle \exp
\rho _{n}\left( \mathfrak{sl}\left( 2,\mathbb{C}\right) \right) \rangle $
such that the group generated by $\Gamma $ is Zariski dense. Then $\Gamma $
generates $\mathrm{Sl}\left( n+1,\mathbb{C}\right) $ as a semigroup.
\end{teorema}

Given the basis 
\begin{equation*}
X=\left( 
\begin{array}{cc}
0 & 1 \\ 
0 & 0%
\end{array}
\right) ,\qquad H=\left( 
\begin{array}{cc}
1 & 0 \\ 
0 & -1%
\end{array}
\right) ,\qquad Y=\left( 
\begin{array}{cc}
0 & 0 \\ 
1 & 0%
\end{array}
\right)
\end{equation*}
of $\mathfrak{sl}\left( 2,\mathbb{C}\right) $ we choose a basis $%
\{v_{0},v_{1},\ldots ,v_{n}\}$ of $\mathbb{C}^{n+1}$ such that 
\begin{equation*}
\rho _{n}\left( H\right) =\mathrm{diag}\{n,n-2,\ldots ,-n+2,-n\},
\end{equation*}
$\rho _{n}\left( Y\right) v_{j}=v_{j+1}$ and $\rho _{n}\left( X\right)
v_{j}=j\left( n-j+1\right) v_{j-1}$. The diagonal matrix $\rho _{n}\left(
H\right) $ is a regular real element of $\mathrm{Sl}\left( n+1,\mathbb{C}%
\right) $ and is contained in the standard Weyl chamber $\mathfrak{a}^{+}$
formed by real diagonal matrices with strict decreasing eigenvalues. Hence $%
G_{1}=\langle \exp \rho _{n}\left( \mathfrak{sl}\left( 2,\mathbb{C}\right)
\right) \rangle $ is in the case covered by Corollary \ref%
{corSemiSimpleGeral}. To apply it we must look at the homotopy properties of
the orbits $G_{1}\cdot b_{\Theta }$ on the minimal flag manifolds of $%
\mathrm{Sl}\left( n+1,\mathbb{C}\right) $, that is on the Grassmannians $%
\mathrm{Gr}_{k}\left( n+1\right) $ of $k$-dimensional subspaces of $\mathbb{C%
}^{n+1}$.

By the choice of the standard Weyl chamber $\mathfrak{a}^{+}$ the origin of $%
\mathrm{Gr}_{k}\left( n+1\right) $ is the subspace $b_{k}$ spanned by $%
\{v_{0},v_{1},\ldots ,v_{k-1}\}$. An easy computation shows that the
parabolic subalgebra $\mathfrak{p}=\mathrm{span}_{\mathbb{C}}\{\rho
_{n}\left( H\right) ,\rho _{n}\left( X\right) \}$ is contained in the
isotropy subalgebra at any $b_{k}$. Since the corresponding parabolic
subgroup $P\subset G_{1}$ is connected ($G_{1}$ is a complex Lie group) it
follows that the isotropy subgroups at $b_{k}$ contain $P$. On the other $Y$
does not belong to the isotropy subalgebras of $b_{k}$. This implies that
the isotropy subgroup at $b_{k}$ of the action of $G_{1}$ is the parabolic
subgroup $P$. Hence all the orbits $G_{1}\cdot b_{k}$ are $G_{1}/P\approx
S^{2}$.

In the sequel it will be proved that the orbits $G_{1}\cdot b_{k}$ are not
contractible in $\mathrm{Gr}_{k}\left( n+1\right) $. This will be done by
showing that a canonical line bundle over $\mathrm{Gr}_{k}\left( n+1\right) $
is not trivial when restricted (pull backed) to $S^{2}\approx G_{1}/P$. This
approach is based on the following proposition about the restriction of the
tautological line bundle $\mathbb{C}^{n}\setminus \{0\}\rightarrow \mathbb{P}%
^{n}$ to the projective orbit $S^{2}\approx \langle \exp \rho _{n}\left( 
\mathfrak{sl}\left( 2,\mathbb{C}\right) \right) \rangle \left[ v_{0}\right] $
in the irreducible representation.

\begin{proposicao}
\label{propfiblinhas2}Let $\rho _{n}:\mathfrak{sl}\left( 2,\mathbb{C}\right)
\rightarrow \mathfrak{sl}\left( n+1,\mathbb{C}\right) $ be the irreducible
representation of $\mathfrak{sl}\left( 2,\mathbb{C}\right) $ on $\mathbb{C}
^{n+1}$, $n\geq 1$. With the notation as above let $G_{1}\cdot \left[ v_{0} %
\right] \subset \mathbb{P}^{n}$ be the projective orbit of the highest
weight space $\left[ v_{0}\right] $. Then $G_{1}\cdot \left[ v_{0}\right] $
is diffeomorphic to $S^{2}$.

Let $\pi :\mathbb{C}^{n+1}\setminus \{0\}\rightarrow \mathbb{P}^{n}$ the
tautological line bundle over $\mathbb{P}^{n}$ and denote by $\tau _{n}$ its
restriction (pull back) to $S^{2}=G_{1}\cdot \left[ v_{0}\right] $. Then $%
\tau _{n}$ is not a trivial line bundle. In fact $\tau _{n}$ is represented
by the homotopy class $n\in \mathbb{Z}=\pi _{1}\left( \mathbb{C}_{\times
}\right) $.
\end{proposicao}

\begin{profe}
For the proof we view $\tau _{n}$ as a bundle $\xi _{1}\cup _{a}\xi _{2}$
obtained by clutching along the equator $S^{1}$ trivial bundles $\xi _{1}$
and $\xi _{2}$ on the north and south hemispheres. Such a bundle is trivial
if and only if the clutching function $a:S^{1}\rightarrow \mathbb{C}_{\times
}$ is homotopic to a point (see Husemoller \cite{hus}, Chapter 9.7).

The set $\sigma =\{[e^{z\rho _{n}\left( Y\right) }v_{0}]\in S^{2}:z\in 
\mathbb{C}\}$ is an open Bruhat cell in the flag manifold $S^{2}=G_{1}\cdot %
\left[ v_{0}\right] $ of $\mathfrak{sl}\left( 2,\mathbb{C}\right) $ since $%
\mathfrak{n}^{-}=\mathrm{span}_{\mathbb{C}}\{Y\}$ is the nilpotent component
of an Iwasawa decomposition. Hence $\sigma $ is open and dense. Actually, $%
\sigma =S^{2}\setminus \{[v_{n}]\}$ because $[v_{n}]$ is the space of lowest
weight. Furthermore the map 
\begin{equation*}
\phi :z\in \mathbb{C}\longmapsto \lbrack e^{z\rho _{n}\left( Y\right)
}v_{0}]\in S^{2}\setminus \{[v_{n}]\}
\end{equation*}
is a chart. This map yields the following section over $S^{2}\setminus
\{[v_{n}]\}$ of the tautological bundle: 
\begin{equation*}
\chi _{1}:[e^{z\rho _{n}\left( Y\right) }v_{0}]\in S^{2}\setminus
\{[v_{n}]\}\longmapsto e^{z\rho _{n}\left( Y\right) }v_{0}.
\end{equation*}
The same way 
\begin{equation*}
\psi :w\in \mathbb{C}\longmapsto \lbrack e^{w\rho _{n}\left( X\right)
}v_{n}]\in S^{2}\setminus \{[v_{0}]\}
\end{equation*}
is a chart and 
\begin{equation*}
\chi _{2}:[e^{w\rho _{n}\left( X\right) }v_{n}]\in S^{2}\setminus
\{[v_{0}]\}\longmapsto e^{w\rho _{n}\left( X\right) }v_{n}
\end{equation*}
is a section over $S^{2}\setminus \{[v_{0}]\}$.

The change of coordinates $\psi ^{-1}\circ \phi :\mathbb{C}\setminus
\{0\}\rightarrow \mathbb{C}\setminus \{0\}$ between the charts is 
\begin{equation*}
w=\psi ^{-1}\circ \phi \left( z\right) =\dfrac{1}{z}.
\end{equation*}
In fact, we have 
\begin{equation*}
e^{z\rho _{n}\left( Y\right) }v_{0}=\left( 1,z,\frac{z^{2}}{2!},\ldots , 
\frac{z^{n}}{n!}\right) \qquad e^{w\rho _{n}\left( X\right) }v_{n}=\left(
p_{n}\left( w\right) ,\ldots ,p_{2}\left( w\right) ,p_{1}\left( w\right)
,1\right)
\end{equation*}
where $p_{j}\left( w\right) $ is a polynomial with positive coefficients and 
$p_{1}\left( w\right) =nw$. If $z\neq 0$ then $\left( 1,z,\frac{z^{2}}{2!}
,\ldots ,\frac{z^{n}}{n!}\right) $ and 
\begin{equation*}
\left( \frac{n!}{z^{n}},\frac{n!}{z^{n-1}},\frac{n!}{2!z^{n-2}},\ldots , 
\frac{n}{z},1\right)
\end{equation*}
span the same line. Hence, if $z,w\neq 0$ then $\left[ e^{z\rho _{n}\left(
Y\right) }v_{0}\right] =\left[ e^{w\rho _{n}\left( X\right) }v_{n}\right] $
if and only if $n/z=nw$, which means that $w=1/z$ as claimed. It follows
that $|\psi ^{-1}\circ \phi \left( z\right) |=1$ if $|z|=1$ so that $%
S^{1}=\{z:\left\vert z\right\vert =1\}$ is the same equator of $S^{2}$ in
both charts.

Now the clutching function $a:S^{1}\rightarrow \mathbb{C}_{\times }$ is
given by $\chi _{2}\left( x\right) =a\left( x\right) \chi _{1}\left(
x\right) $ with $x$ in the equator $S^{1}$. To get it take $w=\psi
^{-1}\circ \phi \left( z\right) $. Then 
\begin{eqnarray*}
\chi _{2}\left( w\right) &=&e^{w\rho _{n}\left( X\right) }v_{0}=\left( \frac{%
n!}{z^{n}},\ldots ,\frac{n}{z},1\right) \\
&=&\frac{n!}{z^{n}}\chi _{1}\left( z\right) =n!w^{n}\chi _{1}\left( z\right)
.
\end{eqnarray*}
Hence, for $x\in S^{1}$, $a\left( x\right) =n!x^{n}$ which is not homotopic
to a point, showing that the bundle is not trivial.
\end{profe}

This proposition implies that the projective orbit $G_{1}\cdot \left[ v_{0} %
\right] \approx S^{2}$ ($G_{1}=\langle \exp \rho _{n}\left( \mathfrak{sl}%
\left( 2,\mathbb{C}\right) \right) \rangle $) in $\mathbb{P}^{n}$ of the
highest weight space is not contractible and hence not contained in an open
Bruhat cell. In fact, if the orbit were contractible then the restriction of
the tautological bundle on it would be trivial.

The same approach yields also the noncontractibility of the orbits $%
G_{1}\cdot b_{k}$ on the Grassmannians $\mathrm{Gr}_{k}\left( n+1\right) $, $%
1\leq k\leq n-1$, where \linebreak $b_{k}=\mathrm{span}\{v_{0},\ldots
,v_{k-1}\}$. To this end we view the Grassmannian $\mathrm{Gr}_{k}\left(
n+1\right) $ as a subset of the projective space $\mathbb{P}\left( \wedge
^{k}\mathbb{C} ^{n+1}\right) $ of the $k$-fold exterior power of $\mathbb{C}%
^{n+1}$. Namely $\mathrm{Gr}_{k}\left( n+1\right) $ is identified to the
projective $\mathrm{\ Sl}\left( n+1,\mathbb{C}\right) $-orbit of the highest
weight space $\left[ \xi _{k}\right] \in \mathbb{P}\left( \wedge ^{k}\mathbb{%
C}^{n+1}\right) $ where $\xi _{k}=v_{0}\wedge \cdots \wedge v_{k-1}$. Via
this identification $G_{1}\cdot b_{k}$ is identified to the projective orbit
of $G_{1}$ through $\left[ \xi _{k}\right] $.

Consider the representation of $\mathfrak{sl}\left( 2,\mathbb{C}\right) $ on 
$\wedge ^{k}\mathbb{C}^{n+1}$ obtained by composing the representation of $%
\mathrm{Sl}\left( n+1,\mathbb{C}\right) $ with $\rho _{n}$. Denote it by $%
\rho _{n}$ as well. Let 
\begin{equation*}
V_{k}=\mathrm{span}\{\rho _{n}^{j}\left( Y\right) \xi _{k}:j\geq 0\}
\end{equation*}
be the $\mathfrak{sl}\left( 2,\mathbb{C}\right) $-irreducible subspace of $%
\wedge ^{k}\mathbb{C}^{n+1}$ that contains $\xi _{k}$. Since \linebreak $%
\rho _{n}\left( X\right) \xi _{k}=0$ and $\xi _{k}$ is an eigenvector of $%
\rho _{n}\left( H\right) $ with eigenvalue $\lambda $ we have $\dim
V_{k}=\lambda +1$. Now, $\lambda $ is the sum of the $k$-largest eigenvalues
of the matrix $\rho _{n}\left( H\right) $, that is, 
\begin{equation*}
\lambda =\sum_{j=0}^{k-1}\left( n-2j\right) =k\left( n-\left( k-1\right)
\right) >0.
\end{equation*}
Hence $\dim V_{k}>0$ and the $\mathfrak{sl}\left( 2,\mathbb{C}\right) $
representation on $V_{k}$ is the non trivial irreducible representation $%
\rho _{m}$ with $m=k\left( n-k+1\right) $.

Now it is clear that the projective orbit $S^{2}=G_{1}\cdot \left[ \xi _{k}%
\right] $ is contained in the projective space $\mathbb{P}\left(
V_{k}\right) $. By Proposition \ref{propfiblinhas2} we have that the
restriction to $G_{1}\cdot \left[ \xi _{k}\right] $ of the tautological
bundle $V_{k}\setminus \{0\}\rightarrow \mathbb{P}\left( V_{k}\right) $ is
not trivial. But the tautological bundle of $\mathbb{P}\left( V_{k}\right) $
is the restriction to it of the tautological bundle of $\mathbb{P}\left(
\wedge ^{k}\mathbb{C}^{n+1}\right) $. We conclude that the restriction to $%
G_{1}\cdot \left[ \xi _{k}\right] $ of the tautological bundle of $\mathbb{P}%
\left( \wedge ^{k}\mathbb{C}^{n+1}\right) $ is a nontrivial bundle. Hence
the orbit $G_{1}\cdot \left[ \xi _{k}\right] \approx G_{1}\cdot b_{k}\approx
S^{2}$ is not contractible in $\mathbb{P}\left( \wedge ^{k}\mathbb{C}%
^{n+1}\right) $. A fortiori it is not contractible in $\mathrm{Gr}_{k}\left(
n+1\right) $.

Hence the group $G_{1}=\langle \exp \rho _{n}\left( \mathfrak{sl}\left( 2,%
\mathbb{C}\right) \right) \rangle \subset \mathrm{Sl}\left( n+1,\mathbb{C}%
\right) $ falls in the conditions of Proposition \ref{propSemiSimpleGeral}
and Corollary \ref{corSemiSimpleGeral}, concluding the proof of Theorem \ref%
{teoParaRepreSl2}.

\section{Semi-simple subgroups containing regular elements}

We consider here two complex semi-simple Lie algebras $\mathfrak{g}%
_{1}\subset \mathfrak{g}$ such that $\mathfrak{g}_{1}$ contains a regular
real element of $\mathfrak{g}$. That is, if $\mathfrak{g}_{1}=\mathfrak{k}%
_{1}\oplus \mathfrak{a}_{1}\oplus \mathfrak{n}_{1}\subset \mathfrak{g}=%
\mathfrak{k}\oplus \mathfrak{a}\oplus \mathfrak{n}$ are compatible Iwasawa
decompositions then there exists a Weyl chamber $\mathfrak{a}^{+}\subset 
\mathfrak{a}$ such that $\mathfrak{a}_{1}\cap \mathfrak{a}^{+}\neq \emptyset 
$. For this case we can apply essentially the same proof of Theorem \ref%
{teoParaRepreSl2} to get analogous semigroup generators of $G$ by subsets
containing $G_{1}=\langle \exp \mathfrak{g}\rangle $.

\begin{teorema}
\label{teoParaSemiComplex}Let $\mathfrak{g}_{1}\subset \mathfrak{g}$ be
complex semi-simple Lie algebras and suppose that $\mathfrak{g}_{1}$
contains a regular real element of $\mathfrak{g}$. Let $G$ be a connected
Lie group with Lie algebra $\mathfrak{g}$ and put $G_{1}=\langle \exp 
\mathfrak{g}\rangle $. Then a subset $\Gamma \subset G$ generates $G$ as a
semigroup provided the group generated by $\Gamma $ is Zariski dense and $%
G_{1}\subset \Gamma $.
\end{teorema}

Before proving the theorem we note that $G$ has finite center because $%
\mathfrak{g}$ is complex, hence the results on semigroups can be applied to $%
G$. Also we can prove the theorem only for $G$ simply connected. For
otherwise we have the simply connected cover $\pi :\widetilde{G}\rightarrow
G $ and $\Gamma $ generates $G$ if and only if $\pi ^{-1}\left( \Gamma
\right) $ generates $\widetilde{G}$. Hence we assume from now on that $G$ is
simply connected.

The proof of the theorem is based on Proposition \ref{propfiblinhas2}. Fix a
Weyl chamber $\mathfrak{a}^{+}\subset \mathfrak{a}$ such that $\mathfrak{a}%
_{1}\cap \mathfrak{a}^{+}\neq \emptyset $. By Corollary \ref%
{corSemiSimpleGeral} we are required to prove that in any flag manifold $%
\mathbb{F}_{\Theta }$ of $G$ the orbit $G_{1}\cdot b_{\Theta }$ is not
contractible where $b_{\Theta }$ is the origin of $\mathbb{F}_{\Theta }$
defined by means of $\mathfrak{a}^{+}$.

To this end we realize the flag manifolds of $G$ as projective orbits in
irreducible representations of $\mathfrak{g}$. If $\omega $ \ is a dominant
weight of $\mathfrak{g}$ (w.r.t $\mathfrak{a}^{+}$) let $\rho _{\omega }$ be
the irreducible representation of $\mathfrak{g}$ with highest weight $\omega 
$. The space of the representation is denoted by $V_{\omega }$ and the
highest weight space is spanned by $v_{\omega }\in V_{\omega }$.

Since $G$ is simply connected it represents on $V_{\omega }$.

It is known that the orbit $G\cdot \left[ v_{\omega }\right] $ of the
highest weight space in the projective space $\mathbb{P}\left( V_{\omega
}\right) $ is a flag manifold of $G$. Precisely, if $\Sigma $ is the simple
system of roots $\Theta =\{\alpha \in \Sigma :\langle \alpha ,\omega \rangle
=0\}$ then $G\cdot \left[ v_{\omega }\right] $ equals $\mathbb{F}_{\Theta }$
as homogeneous spaces.

Our objective is to prove that $G_{1}\cdot \left[ v_{\omega }\right] $ is
not contractible in $G\cdot \left[ v_{\omega }\right] $. Clearly it is
enough to check that $G_{1}\cdot \left[ v_{\omega }\right] $ is not
contractible in $\mathbb{P}\left( V_{\omega }\right) $. We prove this by
applying Proposition \ref{propfiblinhas2} to a representation of a copy of $%
\mathfrak{sl}\left( 2,\mathbb{C}\right) $ inside $\mathfrak{g}_{1}$.

Let $\Pi _{1}$ be the set of roots of $\left( \mathfrak{g}_{1},\mathfrak{a}%
_{1}\right) $. As before if $\alpha _{1}\in \Pi _{1}$ \ then $H_{\alpha
_{1}}\in \mathfrak{a}_{1}$ is defined by $\alpha _{1}\left( \cdot \right)
=\langle H_{\alpha _{1}},\cdot \rangle $. The subspace $\mathfrak{a}_{1}$ is
spanned by $H_{\alpha _{1}}$, $\alpha _{1}\in \Pi _{1}$. Since $\mathfrak{a}
_{1}\cap \mathfrak{a}^{+}\neq \emptyset $ and $\omega $ is strictly positive
on $\mathfrak{a}^{+}$ it follows that there exists $\alpha _{1}\in \Pi _{1}$
such that $\omega \left( H_{\alpha _{1}}\right) \neq 0$.

We choose $\alpha _{1}$ with $\omega \left( H_{\alpha _{1}}\right) \neq 0$
such that $\alpha _{1}$ is positive for the chosen compatible Iwasawa
decomposition, that is, $\left( \mathfrak{g}_{1}\right) _{\alpha
_{1}}\subset \mathfrak{n}^{+}$. We denote by $\mathfrak{g}_{1}\left( \alpha
_{1}\right) $ the subalgebra of $\mathfrak{g}_{1}$ spanned by the root
spaces $\left( \mathfrak{g}_{1}\right) _{\pm \alpha _{1}}$, which is
isomorphic to $\mathfrak{sl}\left( 2,\mathbb{C}\right) $, and put $%
G_{1}\left( \alpha _{1}\right) =\langle \exp \mathfrak{g}_{1}\left( \alpha
_{1}\right) \rangle \subset G_{1}$.

To get a representation of $\mathfrak{g}_{1}\left( \alpha _{1}\right) $ take
a generator $Y$ of $\left( \mathfrak{g}_{1}\right) _{-\alpha _{1}}$ and let $%
V_{\omega }\left( \alpha _{1}\right) $ be the subspace spanned by $\rho
_{\omega }\left( Y\right) ^{j}v_{\omega }$, $j\geq 0$. Since $\left( 
\mathfrak{g}_{1}\right) _{\alpha _{1}}\subset \mathfrak{n}^{+}$ and $%
v_{\omega }$ is a highest weight vector we have $\rho _{\omega }\left(
\left( \mathfrak{g}_{1}\right) _{\alpha _{1}}\right) v_{\omega }=0$. Also $%
v_{\omega }$ is an eigenvector of $\rho _{\omega }\left( H_{\alpha
_{1}}\right) $ with eigenvalue $\omega \left( H_{\alpha _{1}}\right) \neq 0$
. Hence by the usual construction of the irreducible representations of $%
\mathfrak{sl}\left( 2,\mathbb{C}\right) $ the subspace $V_{\omega }\left(
\alpha _{1}\right) $ is invariant and irreducible by $\rho _{\omega }\left( 
\mathfrak{g}_{1}\left( \alpha _{1}\right) \right) $. Since the eigenvector $%
\omega \left( H_{\alpha _{1}}\right) \neq 0$ we have $\dim V_{\omega }\left(
\alpha _{1}\right) \geq 2$ and we get a nontrivial representation of $%
\mathfrak{sl}\left( 2,\mathbb{C}\right) $ on $V_{\omega }\left( \alpha
_{1}\right) $. In this representation $v_{\omega }$ is a highest weight
vector because $\rho _{\omega }\left( \left( \mathfrak{g}_{1}\right)
_{\alpha _{1}}\right) v_{\omega }=0$. A posteriori $\omega \left( H_{\alpha
_{1}}\right) =\dim V_{\omega }\left( \alpha _{1}\right) -1>0$.

Now we can apply Proposition \ref{propfiblinhas2} to the representation of $%
\mathfrak{sl}\left( 2,\mathbb{C}\right) $ on $V_{\omega }\left( \alpha
_{1}\right) $ to conclude that the restriction to $G_{1}\left( \alpha
_{1}\right) \cdot \left[ v_{\omega }\right] $ of the tautological bundle $%
V_{\omega }\left( \alpha _{1}\right) \setminus \{0\}\rightarrow \mathbb{P}%
\left( V_{\omega }\left( \alpha _{1}\right) \right) $ is not trivial. This
implies, by restricting twice, that the tautological bundle $V_{\omega
}\setminus \{0\}\rightarrow \mathbb{P}\left( V_{\omega }\right) $ restricts
to a nontrivial bundle on the orbit $G_{1}\left( \alpha _{1}\right) \cdot %
\left[ v_{\omega }\right] $. Hence the orbit $G_{1}\left( \alpha _{1}\right)
\cdot \left[ v_{\omega }\right] $ is not contractible on $\mathbb{P}\left(
V_{\omega }\right) $ so that it is not contractible on $\mathbb{F}_{\Theta
}=G\cdot \left[ v_{\omega }\right] $. Finally, the orbit $G_{1}\cdot \left[
v_{\omega }\right] $ is not contractible as well since it contains $%
G_{1}\left( \alpha _{1}\right) \cdot \left[ v_{\omega }\right] $. By
Corollary \ref{corSemiSimpleGeral} this finishes the proof of Theorem \ref%
{teoParaSemiComplex}.

\vspace{12pt}%

\noindent%
\textbf{Example:} The standard realizations of the classical complex simple
Lie algebras of types $B_{l}$, $C_{l}$ and $D_{l}$ as subalgebras of
matrices satisfy the condition of Theorem \ref{teoParaSemiComplex} as
subalgebras of the appropriate $\mathfrak{sl}\left( n,\mathbb{C}\right) $.
In fact, for $C_{l}$ and $D_{l}$ one has Cartan subalgebras of diagonal $%
2l\times 2l$ matrices $\mathrm{diag}\{\Lambda ,-\Lambda \}$ with $\Lambda $
diagonal $l\times l$, while for $B_{l}$ the Cartan subalgebra is given by $%
\mathrm{diag}\{0,\Lambda ,-\Lambda \}$. The $\mathfrak{a}_{1}$ are given by
such matrices with real entries. A quick glance to these diagonal matrices
shows that in any case $\mathfrak{a}_{1}$ contains diagonal matrices with
distinct eigenvalues, and hence regular real elements of $\mathfrak{sl}%
\left( n,\mathbb{C}\right) $, $n=2l$ or $2l+1$.

\end{document}